\documentclass[a4paper, 12pt]{article}
\usepackage{graphicx, color}
\usepackage{graphics}
\usepackage{titling}
\usepackage{multicol}
\begin{document}

\title{History of Leningrad Mathematics\\ in the first half of the 20th century}
\author{Alexander I. Nazarov\thanks{St. Petersburg Dept of Steklov Mathematical Institute, Fontanka 27, St. Petersburg,  191023, Russia, and St. Petersburg State University, Universitetskii pr. 28, 
St.Petersburg, 198504, Russia.}\setcounter{footnote}{6} \ and
Galina I. Sinkevich\thanks{St. Petersburg State University of Architecture and Civil Engineering, 2nd Krasnoarmejskaja ul. 4, St. Petersburg, 190005, Russia.}
}

\maketitle

\begin{abstract}
 The first half of the 20th century in the history of Russian mathematics is striking with a combination of dramaticism, sometimes a tragedy, and outstanding achievements. 
 The paper is devoted to St. Petersburg-Leningrad Mathematical School. It is based on a chapter in the multi-author monograph \cite{MatPet}. 
 
\end{abstract}

\section*{Introduction}

 Russian sciences, including mathematics, came into being in St. Petersburg. Peter the Great founded St. Petersburg in 1703. Petersburg Academy of Sciences was created in 1725 
 in accordance with Peter's Decree. Whereas there were no local scientists as yet, academicians were invited from other countries and obligated to transfer their knowledge to 
 Russian students. Mathematicians, mechanicians and astronomers, like L. Euler, brothers D. Bernoulli, and N. Bernoulli, J.N. Delisle, and A.J. Lexell, were among Russia's first academicians.

 In the 19th century, higher learning institutions began to be massively created. Such prominent mathematicians as M.V. Ostrogradsky\footnote{Mikhail V. Ostrogradsky (1801--1861) was an 
 outstanding Russian mathematician, mechanician and organizer of mathematical education, member of the Academy. He made the essential contributions in many fields of mathematics and physics.}, 
 V.Ya. Bunyakovsky\footnote{Victor Ya. Bunyakovsky (1804--1889) was a famous Russian mathematician and mechanician, Vice President of St. Petersburg Academy of Sciences. His key works were devoted 
 to probability theory and its applications, number theory and analysis. Bunyakovsky played an important role in the organization of mathematical education in Russia.}, 
 P.L. Chebyshev\footnote{Pafnuty L. Chebyshev (the right pronunciation is Chebysh\'{o}v) (1821--1894) was an outstanding Russian mathematician and mechanician, founder of the St. Petersburg 
 mathematical school, who is remembered primarily  for his works on number theory, probability and approximation theory.} appeared. 
 In the second half of the 19th century St. Petersburg mathematical school was formed to comprise such areas as the number theory, real and complex analysis, probability theory, 
 differential equations, theoretical mechanics. The names of E.I. Zolotarev\footnote{Egor I. Zolotarev (the right pronunciation is Zolotary\'{o}v) (1847--1878) was a prominent Russian mathematician, 
 student of P.L. Chebyshev, adjunct of the Academy; specialised in the integration theory, complex variable, and number theory; author of one of the simplest demonstrations of the reciprocity 
 law. Died tragically at the age of 32.}, A.A. Markov, A.M. Lyapunov graced the history of mathematics  of this city. Many Russian scientists, 
 who graduated from St. Petersburg University, have maintained the specific features of the Petersburg school in their activities: rigorous algorithmic thinking, engineering 
 set up of the problem, presentation of results in a form convenient for further use.

 The 20th century had leaded to social shocks. Even the name of the city was changed three times: from 1914 to 1924, it was Petrograd; from 1924 to 1991, Leningrad; 
 from 1991 it is again St. Petersburg.
 
 World War I, revolutions, and the civil war disrupted the peaceful  life of this city. However, despite the hardships, fruitful scientific life went on. St. Petersburg may be proud 
 of such names as V.A. Steklov, B.G. Galyorkin, N.M. Gyunter, A.N. Krylov, S.N. Bernstein, V.I. Smirnov, G.M. Fichtenholz, A.A. Friedmann, B.N. Delaunay, N.E. Kochin, P.Ya. Polubarinova-Kochina, 
 D.K. Faddeev, S.G. Mikhlin, S.L. Sobolev, L.V. Kantorovich, A.D. Aleksandrov, Yu.V. Linnik. In defiance of the dramatic events, the success of Leningrad mathematical school suits the names 
 of its founders and forerunners, Euler and Chebyshev.

 \section*{I}

 In the early 20th century, St. Petersburg was the capital of the Russian Empire. The Academy of Sciences, which included mathematical cabinet, was in this city. Petersburg Mathematical School 
 was one of the leading in Europe, the University trained professors for the entire Empire. A strong mathematical community formed in the city: members of the Academy A.A. Markov\footnote{Andrey 
 A. Markov (1856--1922) was an outstanding Russian mathematician, student of P.L. Chebyshev, member of the Academy; author of important results in the number theory, differential equations, 
 function theory and probability theory. He introduced an important class of stochastic processes later named after him.}, 
 A.M. Lyapunov\footnote{Alexander M. Lyapunov (1857--1918) was a great Russian mathematician and mechanician, student of P.L. Chebyshev, member of the Academy. He created the stability theory 
 of the dynamical systems and established several breakthrough results in mathematical physics and probability theory. Lyapunov shot himself dead when his wife died.}, V.A. Steklov\footnote{Vladimir 
 A. Steklov (1864--1926) was a prominent Russian mathematician, mechanician and physicist, Vice President of the Academy of Sciences; originator of Leningrad school of mathematical physics. He holds 
 the credit of preserving the Academy in the post-revolutionary period. He founded the Institute of Physics and Mathematics (1921). Nowadays, the Mathematical Institutes in Moscow and St. Petersburg 
 are named after Steklov.}, and A.N. Krylov\footnote{Alexey N. Krylov (1863--1945) was an outstanding Russian and Soviet mathematician, mechanician and theorist of ship design, member of the Academy of 
 Sciences, three-star navy general. His main areas of research were ship theory, the theory of magnetic and gyro compasses, hydrodynamics, and computational mathematics. He also wrote works devoted 
 to magnetism, gunnery, mathematical physics, astronomy, geodesy, insurance business, history of mathematics, and issues in education. He translated Newton's {\it Philosophiae Naturalis Principia 
 Mathematica} into Russian. Together with V.A. Steklov, Krylov played the most important role in the survival of the Academy of Sciences in the early 1920s.}, professors Yu.V. Sokhotsky\footnote{Julian 
 V. Sokhotsky (1842--1927) was a famous Polish-born Russian mathematician. His works are devoted to the number theory and complex variable theory. He is the author of a theorem on behaviour of analytic 
 function in a neighbourhood of an essentially singular point (1868), and the Sokhotsky--Plemel formulas (1873) which are used in quantum physics to the present time. He was President of St. Petersburg 
 Mathematical Society from 1892.}, N.M. Gyunter\footnote{Nikolai M. Gyunter (1871--1941) was a famous Russian and Soviet mathematician, Corresponding member of the Academy. His main works are devoted 
 to the  theory of differential equations and mathematical physics. Gyunter was the first to provide a rigorous and systematic presentation of the potential theory (1934).}, D.F. 
 Selivanov\footnote{Dmitry F.  Selivanov (1855--1932) was a Russian mathematician and educator, student of Chebyshev and Sokhotsky. After 1917, he was the first democratically elected University 
 rector. However, in 1922, he was  ordered out of Russia for teaching mathematics ``not the way the reds should''. He lived the rest of his life in Prague.}, 
 and A.V. Vasiliev\footnote{Alexander V. Vasiliev (1853--1929) was a Russian mathematician, mathematics historian, and famous public figure. In 1920, he initiated the creation of Petrograd Physical 
 and Mathematical Society. He worked in the history of mathematics, popularization of new mathematical theories in Russia, and organizing scientific life throughout Russia and worldwide.} determined 
 both the research and educational levels,  carrying on traditions of the previous century. In 1890, a Mathematical Society appeared in St. Petersburg (the first President of Society was V.G. 
 Imshenetsky\footnote{Vasily G. Imshenetsky (1832--1892) was a Russian mathematician and mechanician, member of the Academy, expert in partial derivative equations. He was the first President of St. 
 Petersburg Mathematical Society.}, from 1892 it was headed by  Yu.V. Sokhotsky). As of 1913, series of translations, {\it New Ideas in Mathematics}, established by A.V. Vasiliev, were published.
 
% \vskip-0.5cm
\begin{figure}[!htbp]
\centering
\includegraphics[width=0.9\textwidth]{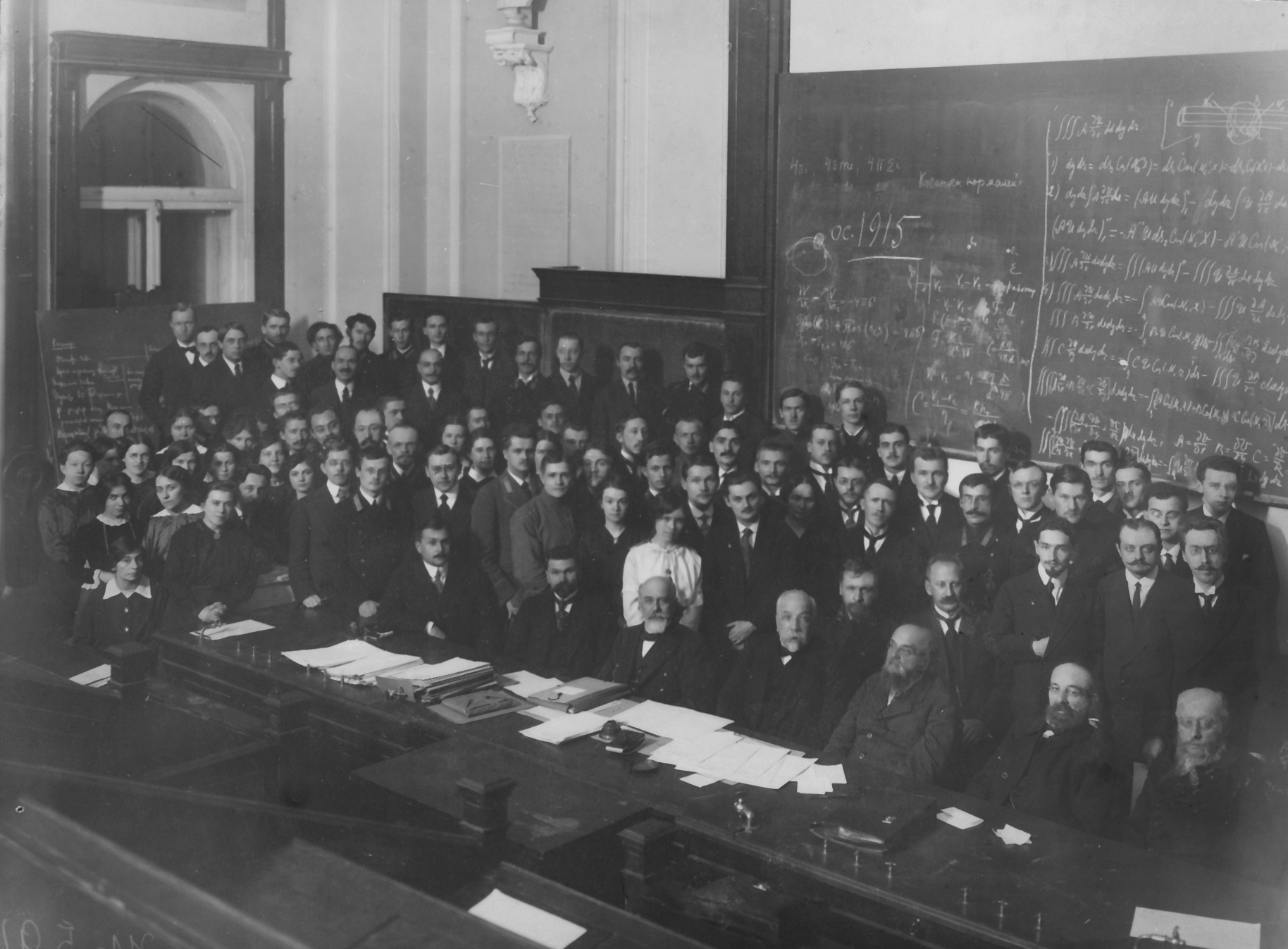}
\caption{\small 1915 year, Imperial Petrograd University. Professors and students of the Mathematical Faculty. Professors E.V. Borisov, H. Ya. Gobi, D.S. Rozhdestvensky, A.I.Voeikov, N.A. Bulgakov, 
Yu.V. Sokhotsky, S.P. von Glazenap, A.A. Ivanov, Kostyakevich.}
\end{figure}

The Europe-wide rise in education in the early 20th century, caused by the need in technical specialists, was conductive to alignment of mathematics taught at secondary schools with 
the progress of science. This was the period of pedagogical activity and high professionalism of teachers of mathematics. The amount of educational and popular scientific literature has grown; 
special magazines appeared in Russia, e.g. {\it Elementary Mathematics Journal, Bulletin of Experimental Physics and Elementary Mathematics}. F. Klein, D. Smith, and E. Borel were leaders 
of this movement in Europe and America; N.Ya. Sonin\footnote{Nikolay Ya. Sonin (1849--1915) was a Russian mathematician, prominent figure in the sphere of organizing education. He worked in St. 
Petersburg as of 1893. His works are devoted to the theory of special functions and its applications.} and K.A. Posse\footnote{Konstantin A. Posse (1847--1928) was a Russian mathematician who played 
a significant role in the reform of mathematical education.} played such part in St. Petersburg. Mathematical circles for pupils, which were run by university teachers, 
appeared at grammar schools. For example, in St. Petersburg Grammar School No.~2, Teachers of Mathematics N.I. Bilibin and Ya.V. Iodynsky established a private mathematical research circle, 
and A.A. Markov guided the steps of some grammar school students in their studies. Having entered the University, the leavers of this school, V.I. Smirnov\footnote{Vladimir I. Smirnov (1887--1974) 
was a famous Russian and Soviet mathematician, member of the Academy; outstanding educator and organiser of science; excellent lecturer, creator of a {\it Course of Higher Mathematics} in five volumes, 
which was reissued numerous times and translated into eight languages. He created an Institute of Mathematics and Mechanics at the University and headed it for 20 years. His main works are devoted 
to the complex variable theory, partial derivative equations, calculus of variations, and to the wave propagation theory. All pre-war Leningrad mathematicians and mathematicians of the first post-war 
years attended his lectures. He also played a great role in the organization of historical mathematical researches in Leningrad.}, 
Ya.D. Tamarkin\footnote{Yakov D. Tamarkin (1888--1945) was a famous Russian and American mathematician. In 1925, he illegally emigrated from USSR, then lived in the United States. Tamarkin's work 
spanned a number of areas of mathematics. He was a proponent and a founding co-editor of the {\it Mathematical Reviews}, together with O. Neugebauer and W. Feller. He was also an active supporter 
of the American Mathematical Society, a member of the council starting 1931, and a Vice President in 1942-43.}, and A.A. Friedmann\footnote{Alexander A. Friedmann (1888--1925) was a prominent Russian 
Soviet mathematician, physicist, and geophysicist, founder of modern physical cosmology. He worked also in atmospheric physics, hydrodynamics and aerodynamics. He discovered non-steady solutions to 
Einstein's equations, which gave rise to the development of the model of non-stationary Universe. Friedmann died of typhoid.}, together with Ya.A. Shochat\footnote{Yakov A. (Yankel) Shohat (1886--1944) 
was a Russian and American mathematician. After 1923, he lived in the USA. Shohat was the first to provide a systematic presentation of general theory of orthogonal polynomials. He was one of the 
publishers of the {\it American Mathematical Society Bulletin} for several years.}, A.S. Bezikovich\footnote{Abram S. Bezikovich (1891--1970) was a prominent Russian and British mathematician. Since 
1924, he illegally emigrated from USSR, then lived in England. His monograph {\it Almost periodic functions} (1932) had won D. Adams prize (Cambridge), and the class of functions he had introduced 
was called Bezikovich functions. He achieved substantial results in the theory of fractal sets.}, and Ya.V. Uspensky\footnote{Yakov V. Uspensky (1883--1947) was a Russian and American mathematician 
and historian of mathematics, member of the Academy. Since 1927, he lived in the USA. His works pertain to the number theory, probability theory, and algebra; research of mechanical quadrature 
convergence; history and methodology of mathematics. He translated J.~Bernoulli's treatise {\it About the Law of Large Numbers} into Russian.}, created a mathematical workshop, where they delivered 
lectures in classical and modern sections of mathematics to each other. 
On the recommendation of V.A. Steklov, they were retained at the University to get ready for professorship\footnote{The getting ready for professorship was similar to the contemporary PhD program.} 
and soon achieved significant research results.

%\vskip-0.5cm
\begin{figure}[!htbp]
\centering
\includegraphics[width=0.8\textwidth]{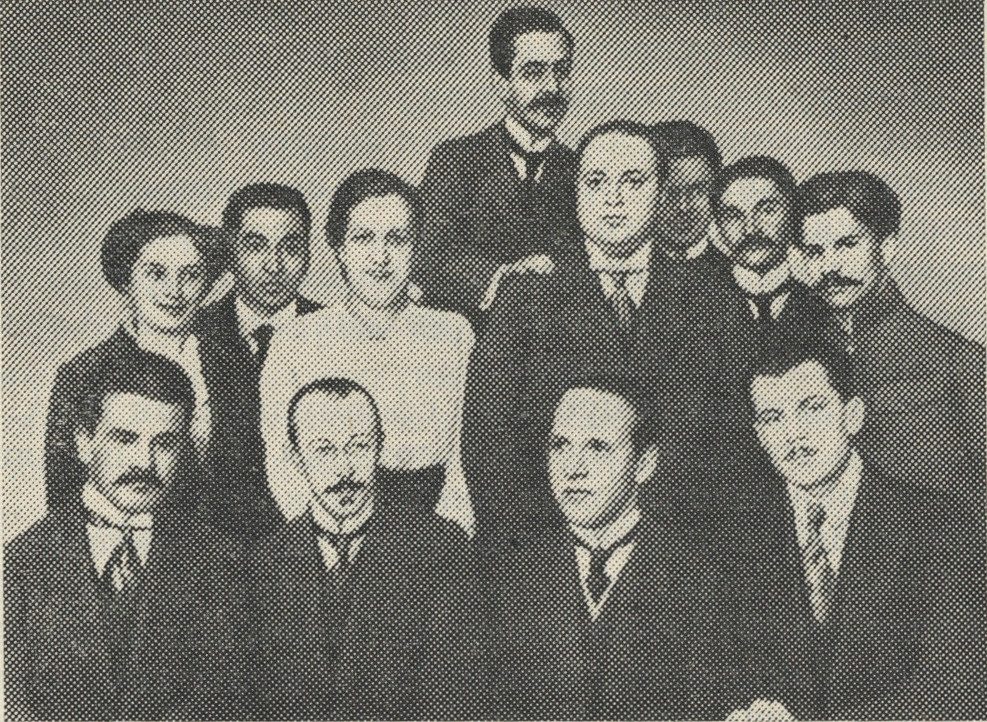}
\caption{\small 1913-1914 years. Left-to-right, first line: unknown (P. Erenfest?), A. Friedmann, G. Weihardt, unknown. Second line: E. Friedmann, A. Shohat, E. Weihardt, Ya. Tamarkin,
unknown, unknown, M. Petelin. Third line: V. Smirnov.}
\end{figure}

The World War I of 1914--1918, the revolution and civil war of 1917--1920 disrupted the untroubled scientific life. The city was flooded by numerous ruined peasants. Authorities assigned 
living quarters to the urban and country poor families in apartments consisting of several rooms, leaving one room to former owners of these apartments. People were coming short of foodstuffs, 
heating, lighting, and bare essentials. Academician Markov complained that he could not attend meetings of the Academy as he had no shoes. The following records about Sokhotsky can be found 
in archives: ``On 10 April 1919, Professor of the First Petrograd University, Yu.V. Sokhotsky, was put on the list of those entitled to a reinforced ration in kind. Rations will be issued at 
the distributing station of the People's Commissariat for Education. Bread will be issued once per week. Please have a bottle for sunflower oil with you.''\footnote{Central State 
Historical Archives of St. Petersburg. F. 14. Schedule 1. No. 6646, 201 leaves. L. 178.} Sokhotsky's letter of 1922: 
``Having spent three winters in succession in an unheated apartment, I have completely upset my health. Now, because of my weakness and in view of the upcoming winter, I have to contact 
KUBU\footnote{KUBU (Commission for the improvement of the life of scientists) was established in 1920 in Petrograd on the writer M. Gorky's initiative. It served as the foundation of the 
first House of Scientists in Russia. There was a boarding house for the elderly scientists.} and ask them to provide me, if it turns possible, with suitable accommodations where I could 
spend winter months in bearable living conditions which pose 
no hazard to life.''\footnote{Central State Historical Archives of St. Petersburg. F. 184. Schedule 2. No. 770, 50 leaves. L. 50.} Classes 
at the University were given in cold rooms without lights. The salary was not paid regularly, often a meager food ration was given instead of money. Most of professors were compelled to teach 
at several learning institutions.

Based on the decrees of the {\it Sovnarkom} (Council of Ministers) passed in 1918--1919, free schools and new rules of enrolment in higher learning institutions were introduced. All working people 
were granted the right to enter any high school regardless of their previous education, without presenting a certificate of secondary education. A bylaw was adopted to grant the right of 
priority access to the universities for workers and peasant poors. In 1920, an illiteracy liquidation campaign began, when primary education was provided to all adults aged 16 to 50. 
In 1919, {\it Rabfak} (the remedial schools for young workers) was formed. In fact, those were preparatory departments designed to provide pre-entry training to those who wished to enter a 
higher learning institution. One could enter this remedial school on recommendation of his trade union or party authorities. New curricula were created to fit the poor level of students' 
training. Academic degrees and titles were abolished (to be reintroduced in 1934).

All these measures caused professors' resistance. For example, Ya.V. Uspensky wrote: ``Taking into account the fact that to succeed at the university, a student should be adequately trained, 
prospective students must be admitted at the university in virtue of their knowledge, not their class affiliation or political commitment.''\footnote{Central State Historical Archives of St. 
Petersburg. F. 7240. Schedule 14. No. 16. L. 185 recto. } In response, the authorities enacted measures 
intended to ``re-educate the bourgeois professors'' and punitive measures like preventive detention, execution, or exile. Only in September and November 1922, 160 scientists were exiled from 
St. Petersburg to Stettin against their will on board the German passenger vessels\footnote{L.D. Trotsky, Chairman of RVSR (Revolutionary Military Council) in 1922, commanded upon this 
campaign: ``We exiled these people as there was no grounds for shooting them dead, but it was impossible to tolerate them.'' In 2003, a memorial sign was installed in St. Petersburg to 
commemorate the place the so called `Philosophers' Ship' departed from.}. Mathematician D.F. Selivanov was among them.

Many scientists emigrated, including mathematicians A.S. Bezikovich, Ya.D. Tamarkin, Ya.A. Shochat, and later on, Ya.V. Uspensky.

\begin{figure}[!htbp]
\centering
\includegraphics[width=0.4\textwidth]{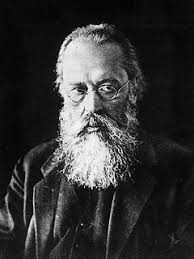}
\caption{\small V.A. Steklov}
\end{figure}
\begin{figure}[!htbp]
\centering
\includegraphics[width=0.4\textwidth]{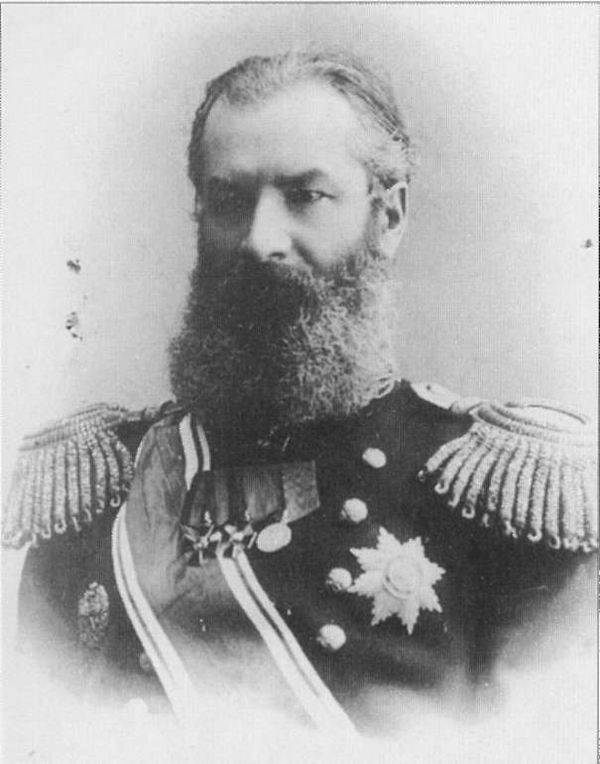}
\caption{\small A.N. Krylov}
\end{figure}

V.A. Steklov, A.N. Krylov, and A.A. Friedmann played an important role in the search for the ways to collaborate with the new authorities and in the preservation of the mathematical community. 
Serving as Vice President of the Academy of Sciences from 1919 to the end of his life, Steklov made an immense contribution in its preservation. In 1921, the Institute of Mathematics and Physics 
was established on his initiative, and he was its first director until 1926 when he died. Thereafter, its directors were, in succession, A.F. Ioffe\footnote{Abram F. Ioffe (1880--1960) was a Russian 
Soviet physicist, outstanding organizer of science, referred to as the ``father of Soviet physics'', member of the Academy, Vice President of the Academy of the USSR (1942--1945). As of 1930s, he 
insisted on the need to conduct intensive nuclear research, which in 1942 marked the start of Soviet nuclear programme. The major merit of A.F. Ioffe was founding a unique physical school which 
enabled the Soviet physics to go global.}, A.N. Krylov, I.M. Vinogradov\footnote{Ivan M. Vinogradov (1891--1983) was a prominent Soviet mathematician, member of the Academy. His works are devoted 
to the analytical number theory. His main achievement was the creation of the trigonometric sums method which is currently one of the key methods in the analytical number theory. With the help of 
this method, he solved some problems which seemed to be out of the reach of the mathematics of the early 20th century (investigation of Waring's problem; solving the ternary Goldbach conjecture 
for all sufficiently large numbers).}.

It was not an easy thing to teach students almost in lack of new textbooks. On assignment of {\it Narkompros} (the Ministry of Education), 
the urgent compilation of new textbooks on mathematics and physics was started. Ya.I. Perelman\footnote{Yakov I. Perelman (1882--1942) was a Russian and Soviet mathematician and physicist, outstanding 
exact sciences promoter and communicator, originator of the recreational science genre. Only in the USSR, his books {\it Physics for Fun, Mathematics for Fun, Astronomy for Fun}, and many others, were 
published 449 times, their total print being more than 13 million copies. In other countries, his books were published 126 times in more than 20 languages. Despite the common misconception, 
Ya.I. Perelman was unrelated to G.Ya. Perelman who proved Poincare's conjecture in 2002--2003.} played an important role in this process. Many of his textbooks were written in the form of 
wonder-books. Ya.I. Perelman also engaged in promotion of the metric weights and measures introduced in Russia in 1918.

The work of the first mathematical society was gradually fading away, and the society ceased to exist after World War I. In 1920, a mathematical circle was created on initiative of A.V. Vasiliev. 
A year later, it was reorganized into Petrograd Physical and Mathematical Society. A.V. Vasiliev was Chairman of the Society for two years; thereafter, it was presided by N.M. Gyunter. 
Among members of this Society were A.S. Bezikovich, S.N. Bernstein\footnote{Sergei N. Bernstein (1880--1968) was an outstanding Russian and Soviet mathematician, member of the Academy. His results and 
methods he had created had great effect on the development of mathematics in the 20th century. He provided the first solution to the 19th Hilbert problem. He laid the groundwork for the constructive 
function theory. When proving Weierstrass theorem, he built polynomials which are now named after him. He was the first to propose axiomatics of probability theory (1917).}, B.G. 
Galyorkin\footnote{Boris G. Galyorkin (1871--1945) was a prominent Russian and Soviet mechanical engineer and mathematician, member of the Academy, engineer-lieutenant-general. He is most famous for 
his results in the elasticity theory. He designed and consulted on large hydroelectric power station projects and projects of other industrial facilities. The approximate method of solving boundary 
value problem for differential equations, he propose in 1915, was named after him.}, V.I. Smirnov, V.A. Steklov, Ya.D. Tamarkin, Ya.V. Uspensky, G.M. Fichtenholz\footnote{Grigory M. Fichtenholz 
(1888--1959) was a Russian and Soviet mathematician; brilliant lecturer, author of a three-volume edition of the {\it Course of Differential and Integral Calculus}. The workshops he had created are underlying 
the Leningrad school of the theory of real variable functions and functional analysis. Most of mathematicians in Leningrad in the pre-war and the first post-war years attended his lectures.}, 
V.A. Fock\footnote{Vladimir A. Fock (1898--1974) was an outstanding Soviet theoretical physicist, member of the Academy. His works pertain to quantum mechanics, quantum electrodynamics, quantum field 
theory, statistical physics, relativity theory, gravitation theory, radio-physics. He also worked on philosophical problems of physics.}, and 
A.A. Friedmann. As of 1926, {\it The Journal of Leningrad Physical and Mathematical Society} founded by Steklov began to be published (after Steklov's death in 1926, the editorial staff of 
this Journal was headed by Ya.V. Uspensky). It was one of the few mathematical journals published in the USSR at that time.

In the middle 1920s, the generation of students enrolled in higher learning institutions after the revolution, completed their education and entered their career, truly believing to the new 
authorities. The new mathematical community was painfully formed in confrontation of the old versus new. In 1927, the All-Russian Congress of Mathematicians was held in Moscow. Soviet 
mathematics was declared to be a `party and class science'. The authorities began chasing `wreckers' in all spheres of life; politically motivated trials began. Old professors vindicated 
their traditions, ethics, the right to freedom of research, the right to freedom of conscience. This struggle reached its climax in 1929. The Ministry of Education issued instructions: 
in spite of the university charter, professors might be dismissed from their positions by resolution of administration. In Moscow, D.F. Egorov\footnote{Dmitry F. Egorov (1869--1931) was 
a famous Russian and Soviet mathematician, worked in Moscow. As of 1910, Egorov conducted mathematical research workshops which paved the way for Moscow school of the theory of real variable 
functions. Since 1921, he was Vice President, and as of 1923, President of Moscow Mathematical Society. In 1929, was persecuted for religious views, and in October 1930, arrested. Together 
with a famous philosopher A.F. Losev, he was targeted in the investigation of the All-Union Counter-revolutionary Organization, Genuine Orthodox Church (catacomb church). He was exiled to Kazan 
for five years even before the proceedings were over. He died in a hospital after he had launched a hunger strike in prison.} was dismissed from his offices in 1929, 
and in 1930, he was arrested.

The authorities set the course for industrialisation. To create strong industry, the country needed skilled engineers. To replace the old `bourgeois' specialists and train the new ones, 
the authorities decided to enrol in higher learning institutions workers sent by the Communist party; facilitate curricula as much as possible; reduce the duration of study and attendance 
at higher learning institutions and technical schools; and give top priority to the ideology component of learning. In opposition to the Academy of Sciences, a Communist Academy was 
established\footnote{After the Academy of Sciences was transferred to Moscow, the Communist Academy was liquidated in 1936  as superfluous.}. Under this Academy, various societies 
of scientists were formed, including the Society of Marxist Mathematicians. Being desperately inferior to the traditional mathematical 
community professionally, Marxist mathematicians converted competition into ideological struggle with `bourgeois' professors, which caused new political trials. Historical, mathematical, 
and methodological issues associated with the new educational programmes formed the battlefield. Textbooks were abolished in schools to be replaced by books of problems (out of the public 
eye, old teachers continued using the classical textbooks by A.P. Kiselev\footnote{Andrei P. Kiselev (1852--1940) was an outstanding author of school textbooks in mathematics. Graduated from 
St. Petersburg University (1875), where he attended lectures of P.L. Chebyshev. He worked at grammar schools and non-classical secondary schools. By the beginning of the 20th century, he had 
created a brilliant line of textbooks in mathematics for grammar schools and non-classical secondary schools. His study books covered practically all school mathematical subjects: arithmetic, 
algebra, geometry, fundamentals of analysis. For more than 60 years they were the most sustainable textbooks for the national school thanks to their high theoretical level, consistency, 
clarity, and brevity.}). `Laboratory team training'\footnote{Each class was divided into teams; each team worked on an assignment; 
the team leader, or at least one team member, reported the progress to the teacher, and the entire team had a pass. This system was abolished in the USSR in 1932.} was introduced, 
and other experiments were carried on. It was only in 1936 that 
scientists managed to protect schooling against incompetent reforms thanks to the presentation G.M. Fichtenholz made at the session of the mathematical group of the Academy of Sciences.

The elections to the Academy of 1929 were protectionistic. The government insisted on electing their nominees, members of the Academy opposed. There was a threat that the Academy will 
be wiped out. A.N. Krylov saved the situation, having insisted that they should comply with the demand of the authorities.\footnote{Krylov said: ``There is nothing to speculate, we must 
do what the government demands us to do $<\dots>$. We have to elect them. Otherwise, the government will send the Academy to the devil together with all academicians.'' \cite[p. 444]{Krylov}.}

In 1930, OGPU\footnote{OGPU (the Joint State Political Directorate) was the name of the secret police of the Soviet Union from 1923 to 1934.} conducted mass arrests of 
``underground anti-Soviet parties'' on fabricated evidence. In Moscow, OGPU ascribed leadership of the `Party of Russia's Renaissance' to academicians 
N.N. Luzin\footnote{Nikolai N. Luzin (1883--1950) was an outstanding Russian and Soviet mathematician, member of the Academy; one of the originators of descriptive theory of sets 
and functions. Together with D.F. Egorov, he founded the Moscow mathematical school and the famous community of Luzin’s students, Luzitania. Luzin’s educational outcome was immense -- it was 
the one-off event when an outstanding scientist nurtured equally outstanding scientists. Many of them created their own scientific schools. Luzin’s life was spoiled by political harassment 
orchestrated around him in 1936. See \cite{Luzin}.} and S.A. Chaplygin\footnote{Sergey A. Chaplygin (1869--1942) was a prominent Russian and Soviet mechanician and mathematician, member of the Academy; 
one of the founders of modern aeromechanics and aerodynamics. S.A. Chaplygin’s key works were devoted to aerohydrodynamics, non-holonomic mechanics, theory of differential equations, and theory 
of aviation.}. These allegations were much-publicized, there was much discussion at meetings and rallies, however, matters didn't come to arrests of the members 
of the Academy.

\begin{figure}[!htbp]
\centering
\includegraphics[width=0.4\textwidth]{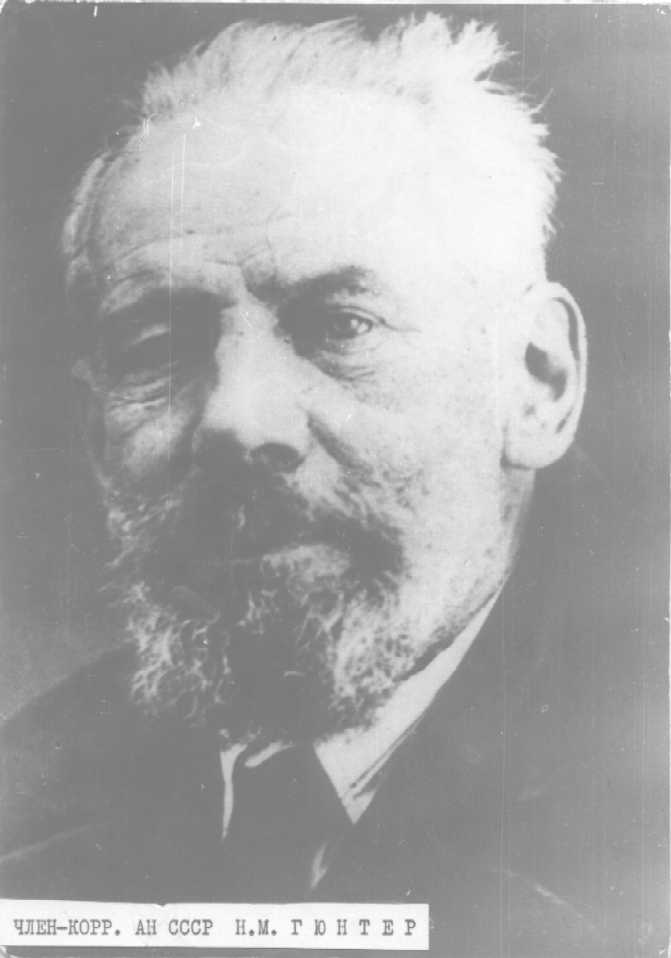}
\caption{\small N.M. Gyunter}
\end{figure}

In Leningrad, professor N.M. Gyunter was the leader of the old school mathematicians. He was the chairman of the Physical and Mathematical Society, Corresponding member of the Academy. 
The Society of Marxist Mathematicians was in opposition to him. They accused the Academy, University, and mathematical community of Leningrad at large, of being separated from practice, 
opposing the teaching reform, `kiselevschina'\footnote{Literally `kiselev-ism'. They meant the use of A.P. Kiselev's textbooks, which were abandoned in the course of the educational reform of 
1918--1932. This reform proved to be unsound. The new curriculum in mathematics of 1935 implied the return to Kiselev's textbooks. } and cliquishness. In 1931, the Society of Marxist 
Mathematicians published a pamphlet entitled {\it ``On the Leningrad Mathematical Front''} \cite{LenMatFront}. 
They faulted `Gyunter's faction' for contemptuous disregard of professionalism of Marxist mathematicians. S.A. Bogomolov (1877--1965), who defended in his works G. Cantor's assertion that 
``the essence of mathematics is in its freedom'' and the right of the mathematician to designate his own scope of research, was also criticised. A.V. Vasiliev, Ya.V. Uspensky, and S.A. Bogomolov 
were criticised for political apathy and academicism of their views on history and teaching of mathematics.

After the hate campaign against Gyunter, on the advice of V.I. Smirnov, the Physical and Mathematical Society was dissolved to save its members from political repressions.

\begin{figure}[!htbp]
\centering
\includegraphics[width=0.7\textwidth]{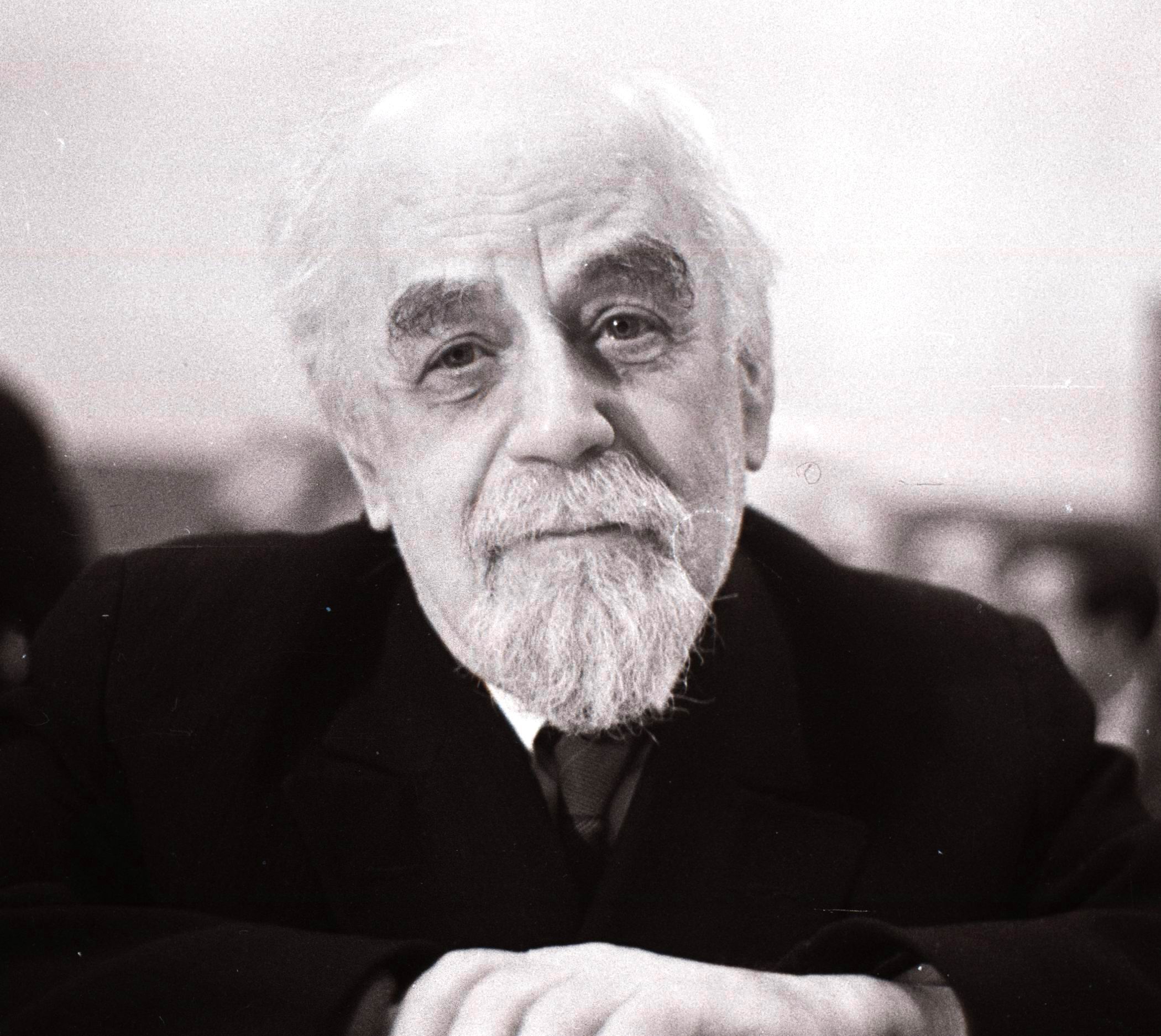}
\caption{\small V.I. Smirnov}
\end{figure}

Gyunter was forced to write a letter of contrition.\footnote{It should be noted that the ``red professor'' L.A. Leifert, main initiator of the hate campaign against Gyunter, was later, in 1938, 
arrested in Voronezh, without due grounds accused of ``participating in a terrorist sabotage and wrecking organisation'' and shot dead.} Having given up the chair, he managed to escape arrest. 
However, he remained at the University as a professor and could continue his scientific activity. Gyunter died on 4 May 1941, having left his mathematical library for the Institute of Mathematics.

In the end of 1933, the Academy of Sciences, which was formerly within the jurisdiction of the CIK\footnote{CIK (the Central Executive Committee of the Soviet Union) was the highest governing body 
in the Soviet Union in the interim of the sessions of the Congress of Soviets of the Soviet Union, existed from 1922 until 1938.}, was transferred under the jurisdiction of the 
Sovnarkom\footnote{Government of the USSR in the period from 1923 to 1946.}. It was resolved to move its 
institutions, including the Institute of Physics and Mathematics, from Leningrad to Moscow in order to centralize research and publications, to make it closer to government, and to facilitate 
science management. The library of the Academy of Sciences, Archives, Publishing House, and some institutes affiliated with the Academy, were left in Leningrad. On 28 April 1934, Presidium of 
the Academy resolved to divide the Institute of Physics and Mathematics into the Institute of Physics and Institute of Mathematics. I.M. Vinogradov was elected as a director of the latter one.

In 1936, the so-called ``Pulkovo Case''\footnote{``Pulkovo Case'' was a criminal case commenced against a group of Soviet scientists based on a baseless accusation of ``participating in a terrorist 
organisation whose goal was to overthrow the Soviet government and establish fascist dictatorship in the USSR''. The first large group of arrests in autumn 1936 affected astronomers of Pulkovo 
Observatory. Therefore, the case was called so. 12 people were shot dead, including B.P. Gerasimovich, director of Pulkovo Observatory, and B.V. Numerov, professor of 
the University, Corresponding member of the Academy, founder and director of several astronomic institutes.} was framed up in Leningrad. As a result, more than one hundred employees of 
Leningrad astronomic institutions were arrested and sentenced to various terms of imprisonment.

Political trials affected all segments of Soviet society. However, the scientific and engineering and technical community, as well as party and top military elite, was in the first place in these 
trials. Investigations were initiated by colleagues' and students' delations, declarations in newspapers and at meetings; thereafter, OGPU pursued the case. Using heavy-handed approach and 
provocative investigative interview methods, they compelled the alleged offender not only to confess to untruths, but to bear false witness against the ambit of acquaintances; thus the number 
of alleged offenders expanded.

\section*{II}

In spite of the stressful environment, mathematicians continued working hard.

A.V. Vasiliev wrote a historical sketch on integer numbers (1919) and a book entitled {\it ``Mathematics. Issue 1 (1725-1826-1863)''} (1921), the first significant work on the development of 
mathematics in Russia, where he analysed in detail works of Lobachevsky\footnote{Nikolai I. Lobachevsky (1792--1856) was an outstanding Russian mathematician, one of the originators of non-Euclidean 
geometry; rector of Kazan University (1827--1846).} and Chebyshev. Initiated by Vasiliev and Steklov, a formal meeting was held in 1921 at the Academy devoted to the 100th anniversary 
of Chebyshev.

In 1922-24, A.A. Friedmann's results appeared based on his research of relativistic models of the universe and of non-stationary solutions of Einstein's equations. In 1922, Friedmann published 
an article entitled {\it ``On the curvature of space''} and his lithographed thesis, {\it ``Essay on Compressible Fluid Mechanics''}, which laid the groundwork for theoretical meteorology. 
In 1923, his famous book {\it ``The World as Space and Time''} was published.

G.M. Fichtenholz worked in St. Petersburg as of 1915; S.N. Bernstein, as of 1933. V.I. Smirnov, R.O. Kuzmin\footnote{Rodion O. Kuzmin (1891--1949) was a Soviet mathematician, Corresponding member 
of the Academy. His main works are devoted to the number theory, analysis, probability theory, and elasticity theory. Together with N.M. Gyunter, they were co-authors of the famous {\it Book of 
Problems in Higher Mathematics} (1930s).}, A.A. Markov Jr.\footnote{Andrey A. Markov Jr. (1903--1979) was a famous Soviet mathematician, Corresponding member of the Academy; founder of the Soviet 
school of constructive mathematics; son of outstanding Russian mathematician A.A. Markov.}, and physicist V.A. Fock played important roles in the 
mathematical community. A new generation came in 1920s: I.A. Lappo-Danilevskij\footnote{Ivan A. Lappo-Danilevsky (1895--1931) was a Soviet mathematician, student of V.I. Smirnov, Corresponding member 
of the Academy. His key works are devoted to the theory of analytic matrix functions and to its applications in the theory of linear differential equations, where he achieved a number of fundamental 
results. The scientist died of acute heart failure during a professional business trip in Germany.}, N.E. Kochin\footnote{Nikolai E. Kochin (1901--1944) was a famous Soviet mathematician and physicist, 
member of the Academy, one of the creators of modern theoretical meteorology. He died in Moscow of a serious disease.} and P.Ya. Polubarinova-Kochina\footnote{Pelagea Ya. Polubarinova-Kochina 
(1899--1999) was a famous Soviet expert in fluid mechanics, member of the Academy; spouse of N.E. Kochin. Her works are devoted to continuum mechanics, hydrodynamics, history of mathematics. As of 
1935, she worked in Moscow and later, in Novosibirsk.}, G.M. Goluzin\footnote{Gennady M. Goluzin (1906--1952) was a Soviet mathematician specializing in the complex variable theory; author of the 
widely known monograph {\it Geometric Complex Variable Theory}.}, I.P. Natanson\footnote{Isidor P. Natanson (1906--1964) was a famous Soviet mathematician, student of G.M. Fichtenholz; founder of 
Leningrad school in constructive function theory; brilliant lecturer. His son, Garald I. Natanson (1930--2003), also was a mathematician specialized in constructive function theory.}, 
S.G. Mikhlin\footnote{Solomon G. Mikhlin (1908--1990) was a famous Soviet mathematician working in analysis, integral equations and computational mathematics. He is best known for the introduction 
of the concept of symbol of a singular integral operator.}, S.L. Sobolev\footnote{Sergey L. Sobolev (1908--1989) was an outstanding Soviet mathematician who paved the way for a number of new research 
areas; member of the Academy; student of N.M. Gyunter and V.I. Smirnov. Based on the notion of distributional derivative, he introduced new function spaces which were later called Sobolev spaces. 
He also developed the theory of distributions. As of 1934, he worked in Moscow. He was among the founders of the Siberian Department of the Academy.}, D.K. Faddeev\footnote{Dmitry K. Faddeev 
(1907--1989) was a prominent Soviet mathematician, Corresponding member of the Academy, leader of Leningrad school of algebra; creator of the theory of cohomology groups, who substantially 
developed Galois theory; expert in numerical analysis. His wife, Vera N. Faddeeva (Zamyatina, 1906--1983), was a famous specialist in numerical linear algebra. His son, Ludvig D. Faddeev (1934--2017) 
was an outstanding Soviet and Russian mathematician, member of the Academy, founder of Leningrad school of modern mathematical physics, President of International Mathematical Union (1987--1990).}, 
L.V. Kantorovich\footnote{Leonid V. Kantorovich (1912--1986) was an outstanding Soviet mathematician, member of the Academy; author of fundamental works devoted to functional analysis and numerical 
analysis; one of the originators of linear programming and its applications in economics; Nobel laureate in economics (1975). As of 1960, he worked in Novosibirsk; as of 1971, in Moscow.}; 
and in the late 1930s -- early 40s, A.D. Aleksandrov\footnote{Alexander D. Aleksandrov (1912--1999) was one of the greatest geometers in the 20th century, member of the Academy; founder of 
Leningrad school of geometry in the large; author of outstanding achievements in geometry, theory of partial differential equations, and mathematical crystallography; Rector of Leningrad University 
(1952--1964); from 1964 to 1986, worked in Novosibirsk.}, S.M. Lozinsky\footnote{Sergey M. Lozinsky (1914--1985) was a famous Soviet mathematician. His research was mainly focused on two areas -- 
constructive function theory and numerical analysis.}, Yu.V. Linnik\footnote{Yury V. Linnik (1915--1972) was an outstanding Soviet mathematician, member of the Academy; author of important achievements 
in the number theory, probability theory, and mathematical statistics. He was a founder of Leningrad school of Probability and Statistics.}, B.Z. Vulikh\footnote{Boris Z. Vulikh (1913--1978) was 
a Soviet mathematician specializing in functional analysis; student of G.M. Fichtenholz; author of famous textbooks in the theory of functions of real variable and functional analysis. His father 
and grandfather were famous experts in mathematical education.}.

\begin{figure}[!htbp]
\centering
\includegraphics[width=0.7\textwidth]{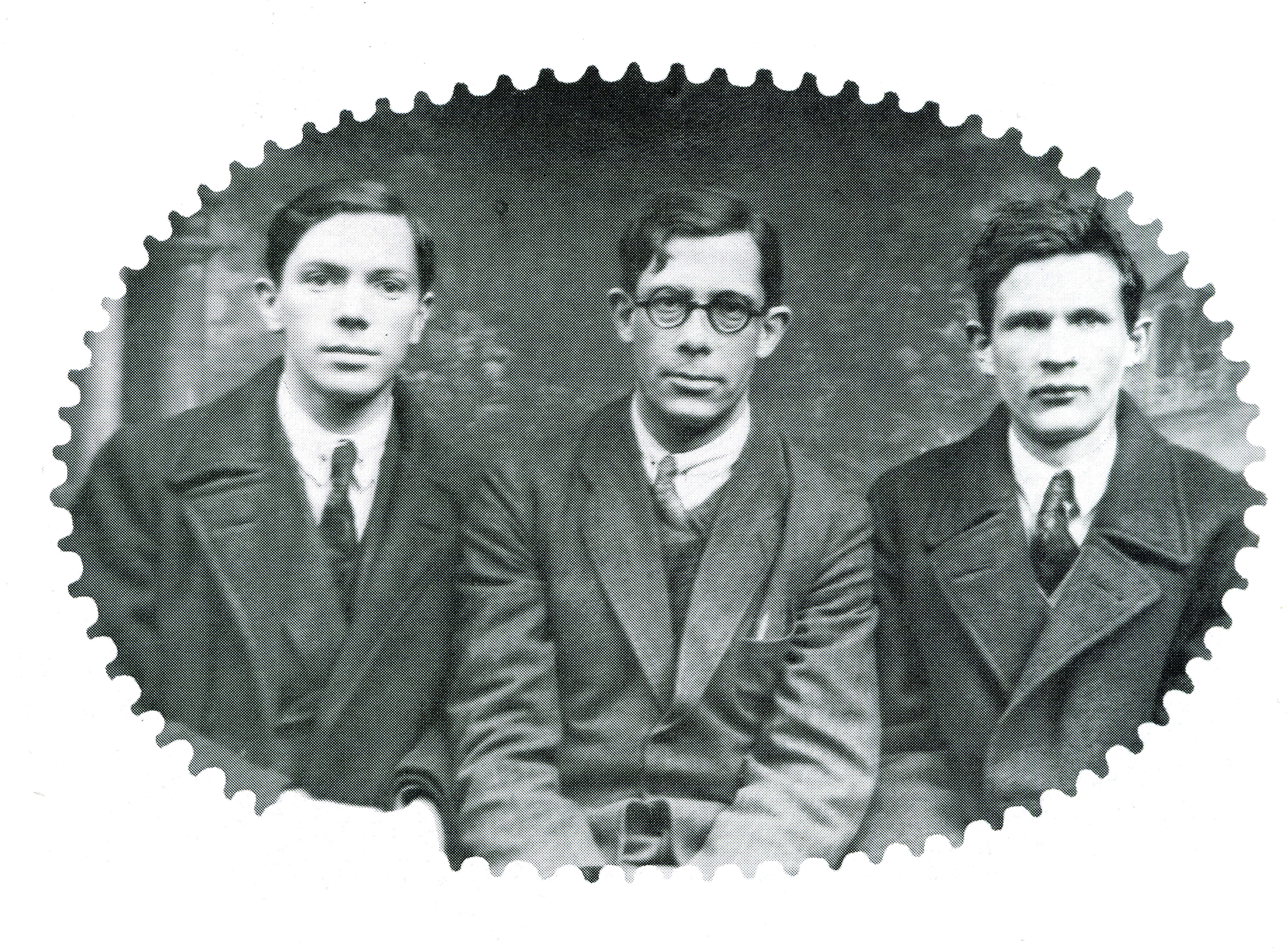}
\caption{\small L.V. Kantorovich, I.P. Natanson and D.K. Faddeev, 1938}
\end{figure}

Before the World War II, many branches of mathematics and applications were actively developed in Leningrad. We can list complex variable theory, theory of functions of real variable, constructive 
function theory, functional analysis, ordinary and partial differential equations, numerical analysis, number theory, algebra, geometry, probability theory, as well as the classical mechanics, 
fluid and gas dynamics, theoretical meteorology, astronomy, ballistics, and quantum mechanics.

In 1934, the Second All-Union Congress of Mathematicians was held in Leningrad. In 1940, Leningrad Department of Steklov Mathematical Institute was set up.

After 1933, many West European scientists came to Russia. Those were communists and communist-leaning, those were also Jewish intimidated by fascist regime. In particular, C. Burstyn\footnote{Celestin 
Burstin (1888--1938) was born in Ternopol, graduated from Viennese University, worked in Austro-Hungary. In 1931, he worked at Leningrad State University and thereafter, he was Director of the 
Institute of Physics and Technology in Minsk. In 1937, he was arrested and died in the prison hospital.} 
(differential equations, differential geometry, algebra), H. M\"{u}ntz\footnote{Herman (Haim) Müntz (1884--1956) was born in Lodz, graduated from Berlin University, was teaching in Germany and 
Poland, worked together with A. Einstein. As of 1929, he was a professor at Leningrad University although remained to be a citizen of Germany. On 1937, he was exiled from the USSR and lived in Sweden.} 
(differential equations), and S. Cohn-Vossen\footnote{Stephan Cohn-Vossen (1902--1936) was born in Breslau (Wroclaw), where he graduated from the university; as of 1930, was teaching at the University 
of Cologne and as of 1933, in Zurich. With D. Hilbert, he was an author of {\it Geometry and the Imagination}. In 1934, he moved to the USSR. He worked in Leningrad State University and in Mathematical 
Institute. He died in Moscow of pneumonia.} (differential geometry) worked for Leningrad University.

As of 1934, initiated by B.N. Delaunay\footnote{Boris N. Delaunay (Delone, 1890--1980) was a famous Soviet mathematician, Corresponding member of the Academy; worked in the area of algebra, number 
theory, computational geometry, mathematical crystallography, and history of mathematics. Since 1934 lived in Moscow.}, mathematical competitions began to be held at schools. G.M. Fichtenholz arranged 
a mathematical circle for school children at the University; leading mathematicians gave lectures for them. In 1934 Ya.I. Perelman, together with colleagues, founded the Russia's first interactive 
museum for children -- The House of Recreational Science.\footnote{There were more than 500 major exhibits and numerous smaller ones (film transparencies, mock-up models, instruments, schemes, and 
diagrams). All of them were grouped in four departments (1939): astronomy (together with meteorology), geography (together with geology), mathematics, and physics (with a room of optics). In 1940, 
they opened new departments: Electricity and Jules Verne Hall. In 1941, the museum was closed, its employees went off to war. Almost the entire exposition was destroyed during the Siege of Leningrad.}

\begin{figure}[!htbp]
\centering
\includegraphics[width=0.4\textwidth]{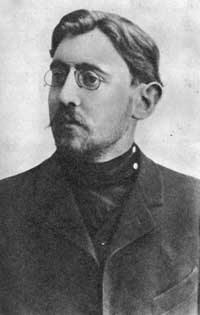}
\caption{\small Ya.I. Perelman}
\end{figure}

As the new requirements to the higher education were introduced in USSR, new courses were developed to combine apprehensibility with sufficient rigor. Many textbooks were written in Leningrad 
before the war, addressed various sections of mathematics and mechanics. Several of them became classical ones, were republished many times and translated into other languages.

In 1938, Presidium of the Academy organized a Commission for the History of the Academy in Leningrad. This Commission was headed by S.I. Vavilov\footnote{Sergey I. Vavilov (1891--1951) was an 
outstanding Soviet physicist, founder of the Soviet school of physical optics, member of the Academy, President of the Academy of the USSR in 1945--1951; younger brother of Nikolay I. Vavilov 
(1887--1943), outstanding biologist who was arrested and died in prison.}. In the first half of 1930s, the translation 
and publishing of mathematical classics was resumed. As a rule, these translations were supplied with detailed comments. The works of Kepler, Galileo, Cavalieri, Descartes, Fermat, Newton, 
L'Hospital, Euler, L. Carnot, Monge, Galois, Dirichlet, etc., were published. In 1932--1938, Leningrad Institute of History of Science and Technology published 10 volumes of its Archives and 
studies in honour of Euler (1935) and Lagrange (1937).

\section*{III}

On 22 June 1941, fascist Germany attacked the Soviet Union. The striking force of German troops was moving towards Leningrad; Finnish troops were approaching from the North. 
On 8 September 1941, they tightened the assault. Leningrad was blocked from all transport directions except for the Ladoga Lake. The siege was overcome but in January 1943 to be finally raised 
on 27 January 1944.

\begin{figure}[!htbp]
\centering
\includegraphics[width=0.8\textwidth]{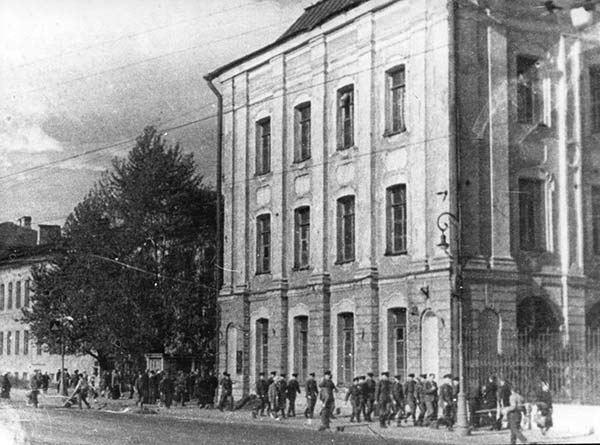}
\caption{\small Leningrad University in 1941}
\end{figure}

Mobilisation began in the very first days of the war. Teachers, students, employees of institutes went up the line. 2,500 people -- employees, postgraduates, students -- were conscripted 
from the University alone. A great number of volunteers signed up in Leningrad people's volunteer army too\footnote{Only on 22 and 23 June 1941, around 100,000 people came to the collecting 
stations of Leningrad Military Registration and Enlistment Office. Seven squadrons were formed of students and teachers of Leningrad State University -- all in all 1,671 men \cite{blokada}.}.

Many people were engaged in building defences around the city. Academician B.G. Galyorkin played the most important role in organizing this process.

\begin{figure}[!htbp]
\centering
\includegraphics[width=0.4\textwidth]{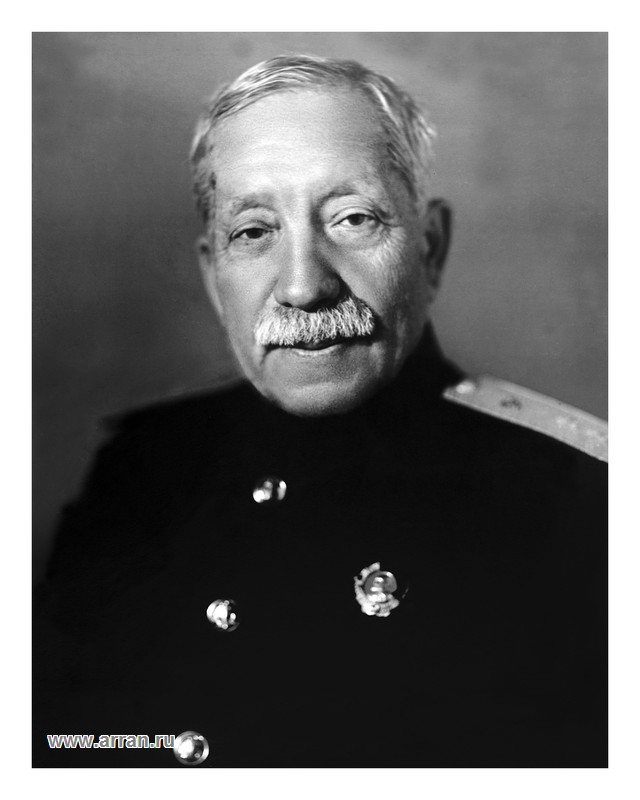}
\caption{\small B.G. Galyorkin}
\end{figure}

By no means all of the scientists' requests to put them into the field were satisfied. In the course of the war, some scientists were recalled from the battlefields and were able to continue 
their research.

Right from the first days of the war, all institutes, museums, libraries, and archives took prompt action to preserve valuable scientific assets -- collections of books and manuscripts, 
archive documents, and instruments. The personnel (mostly women) packed it in boxes and bags, and fetched them to the ground floors and cellars. Fire-safety measures were also taken. 
Manuscripts of M.V. Lomonosov and J. Kepler, and other unique documents from the Archives of the Academy were evacuated to Sverdlovsk (at Ural) together with valuables from the Hermitage. 
In October 1941, under hostile fire, the staff of Astronomical Observatory rescued the instruments and collections of books of the Observatory remaining in Pulkovo.

In the beginning of July 1941, the city's major scientists were taken to the country's safe regions. A.N. Krylov and S.N. Bernstein were evacuated among others. In August 1941, 
the evacuation was suspended to be resumed only after Leningrad was sieged. In autumn 1941, some of the scientists were evacuated by air, and in winter 1942, according to the standard procedure -- 
overice across the Ladoga Lake. The evacuation continued in summer 1942 as well.

The University was substantially sent to Kazan and Yelabuga\footnote{Yelabuga is a small town in the Republic of Tatarstan, located 200 kilometers east from Kazan.}, where, under 
the authority of the Vice-Rector of Leningrad State University V.A. Ambartsumyan\footnote{Victor A. Ambartsumyan (1908--1996) was a prominent Soviet astrophysicist, founder of the Soviet school 
of theoretical astrophysics; member of the Academy; President of the Academy of Armenian SSR (1947--1993); President of International Astronomical Union (1961--1963).}, 
they set up a branch of the University. In March 1942, some faculties were evacuated to Saratov, on Volga, where, under the direction of Dean K.F. Ogorodnikov, classes for students were resumed.

Being evacuated, Leningrad scientists turned their focus toward the needs of the defence industry, making every effort to resume and maintain educational process and intellectual assets for the 
benefit of the victory.

In Yelabuga, astronomer V.V. Sobolev, who later became a member of the Academy, solved problems related to the object visibility calculation, professor V.V. Sharonov drew up tables for 
such calculations. Under the direction of V.I. Smirnov, mechanical engineers (I.P. Ginsburg, M.A. Kovalev, P.G. Makarov) accomplished works in the theory of motion and rotation of 
jet-powered vehicles.

In Saratov, employees of the University conducted researches on propagation of explosive wave in various media (A.A. Grib), stability of anisotropic plates (S.G. Lekhnitsky). In Yaroslavl, 
L.V. Kantorovich performed works on optimal mine disposition and metal cutting.

In 1943, the Academy held a formal meeting in Moscow to celebrate the tercentenary of Newton. Also a book ``Isaac Newton (1643--1943). Collection of Articles to Commemorate His 300th Birthday'' 
was published. A.N. Krylov and other Leningrad scientists were among the authors of this work. In 1943, a Commission for the History of Mathematical Sciences was formed to be presided by A.N. Krylov.

\begin{figure}[!htbp]
\centering
\includegraphics[width=0.8\textwidth]{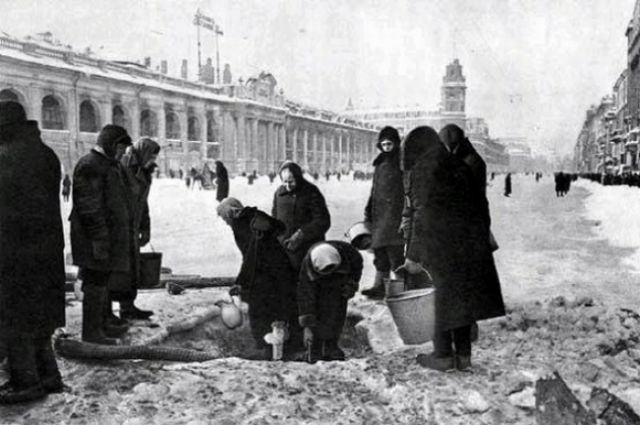}
\caption{\small Cold winter of 1941 in Leningrad. Nevski prospect, people take water from urban water supply}
\end{figure}

The Siege of Leningrad lasted for almost nine hundred days. Day by day, it was harder and harder to live and work in Leningrad. Bombardments were followed by starvation. As of 13 November 1941, 
subject to presentation of redeemable coupons, 300 grams of bread per day was provided to workers and 150 grams, to office employees including research workers. However, even this bread allowance 
was reduced, and as of 20 November, workers were provided 250 grams per day and office employees, 125 grams. The bread allowance was increased, although bitterly moderately, on 25 December 1941: 
350 grams of bread per day began to be provided to workers and 200 grams, to office employees. As of 24 January 1942, this allowance was increased again: 400 grams of bread per day for workers 
and 300 grams for white collars. 98\% of the population suffered dystrophy.

The winter of 1941/42 was the coldest winter of the 20th century in the city (the air temperature fell to $-34$ degrees Celsius). Houses were not heated, water pipes were frozen. Leningraders 
had to take water from ice holes on the Neva, Fontanka, and other rivers. During the bombardments and artillery attacks, window glasses were broken by bomb blasts. People had to close up 
windows with plywood or just curtain them. Electric power was strictly limited. In December 1941, public transport stopped.

More than four thousand people died every day in the city. There were days when up to seven thousand people died. All in all, around eight hundred thousand citizens were starved or frozen 
to death over the period of the siege, and around seventeen thousand people were killed during air strikes and artillery attacks. Ya.I. Perelman, professors of mathematics B.M. Koyalovich 
and N.N. Gernet, and the founder of the first photoelasticity lab in the USSR L.E. Prokofieva-Mikhailovskaya were among them.

Leningraders helped anti-aircraft defence of the city by organizing fire-safe teams. Thus, students and employees of the library of the Academy of Sciences, who kept the watch of its roof, 
saved the library from fire when it was attacked with fire bombs on the night of 10 September 1941. Not only did the fragile women, who worked in the library, weakened by cold and dystrophy, 
manage to preserve its collections of books, they replenished its funds, having preserved libraries of deceased scientists. The Archives of the Academy protected historical documents; its 
employees delivered documents and materials of deceased scientists to the Archive, sometimes carrying them on their own shoulders or on sledges.

\begin{figure}[!htbp]
\centering
\includegraphics[width=0.8\textwidth]{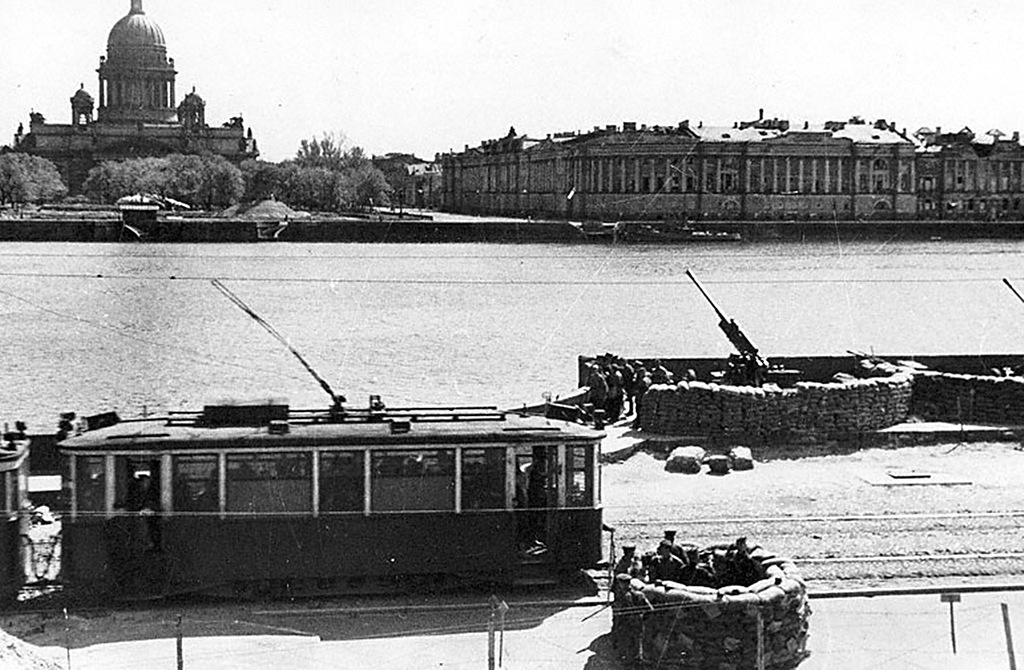}
\caption{\small Anti-aircraft guns on University embankment, 1942}
\end{figure}

On the 1st of September 1941, studies began at the University, first years of education began at the First and Second Medical Institutes, Paediatric Institute, Institute of Civil Engineering, and 
at other higher learning institutions. Senior students began their studies earlier, on the 1st of August. Studying at the institutes and universities, students worked at the same time at 
various enterprises and hospitals. In September--October, 40 higher learning institutions were functioning in the city.

\begin{figure}[!htbp]
\centering
\includegraphics[width=0.8\textwidth]{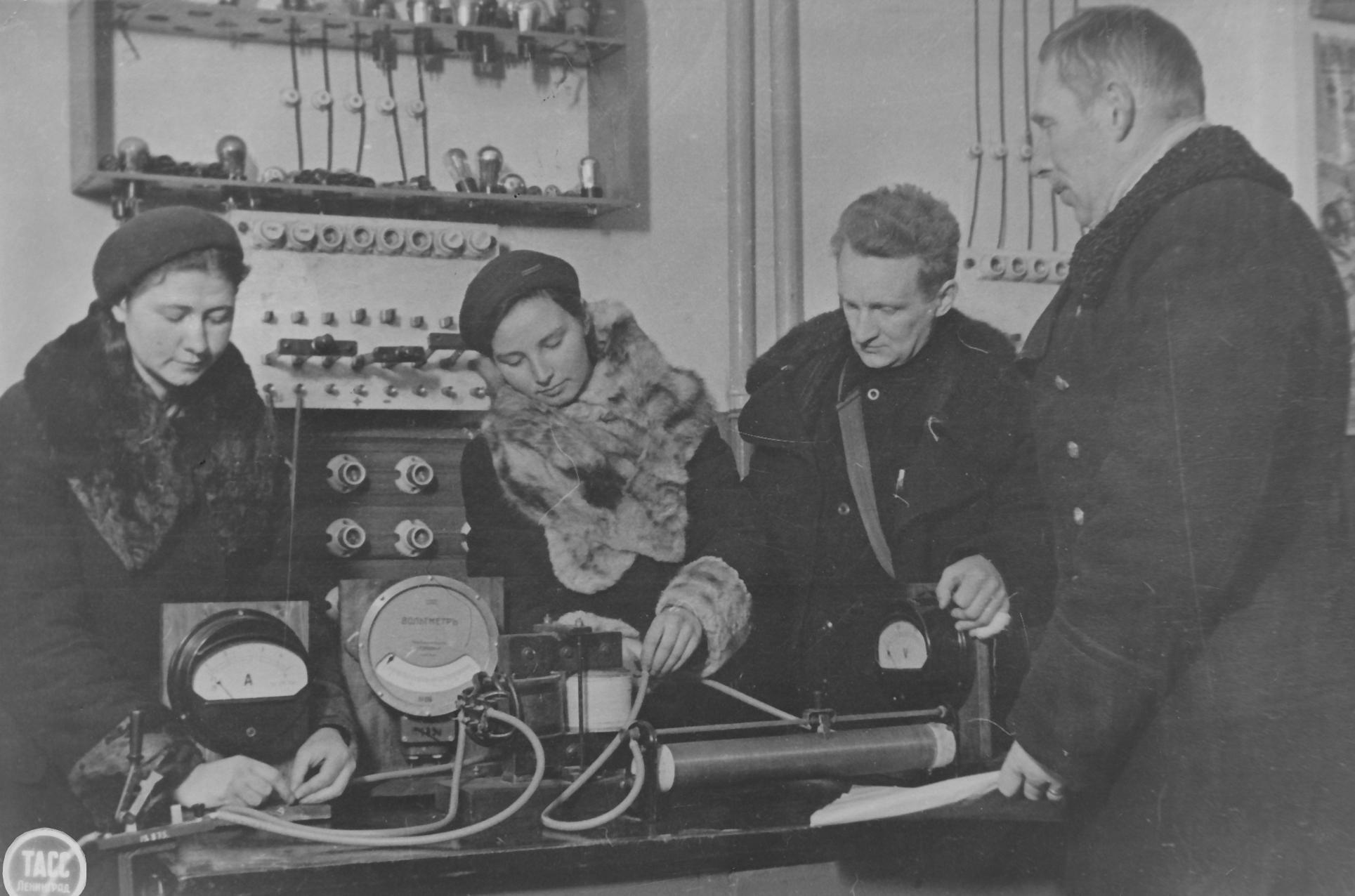}
\caption{\small December of 1941, end-of-term test on experimental physics in Pedagogical Institute}
\end{figure}

Academic and sectoral research institutes did not stop their research work. In July 1941, they established a Defence Proposals Implementation Commission, academician B.G. Galyorkin to 
become its member. This Commission dealt with the challenges related to the creation of effective ways to protect ships against mines, technical support of supplying the city with food 
and evacuating its residents by the Ladoga Lake, and many other issues. Headed by professor I.D. Zhongolovich, a group of employees of the Institute of Astronomy, who remained in the city, 
was engaged in preparing astronomical almanacs and drafting various navigation and ballistic tables. Academician V.A. Fock stayed in Leningrad in the first several months of the siege. 
He and his collaborators were engaged in calculation of ballistic tables and performed various other kinds of defence work. In 1943, professor B.N. Okunev published in the sieged city 
three monographs on ballistics. In 1941--1942, Ya.I. Perelman gave lectures on inland navigation without instruments for reconnaissance men from Leningrad front and Baltic Fleet, and 
for guerilla fighters. The House of Scientists kept working.

L.V. Kantorovich taught algebra, geometry, and probability theory at the Higher Naval Engineering Technical School. In winter, the exploration of the Ladoga Lake ice sheet and investigation 
of opportunities of its use for military purposes were of vital importance for the sieged Leningrad. Food was delivered to the besieged Leningrad by lorries which were moving overice across 
the Ladoga. They used to call this route ``The Road of Life''. A group of scientists including Kantorovich calculated the ice sheet strength under load. They viewed the ice as an engineering 
structure; when calculating the strength, they took into account quality of the material, meteorological situation, and load conditions. In 1943, they derived formulae to calculate the 
load-bearing capacity of the ice sheet in various situations: under stationary load, slowly moving load, movement at an average speed. Essentially, they found the optimal solution for 
movement depending on the changes in the environment. Results of this research found a striking application in practice in December 1943 -- January 1944. Recommendations of L.V. Kantorovich 
and S.S. Golushkevich enabled the vehicles to force through the ice at maximum speed beyond its sturdiness: propeller-driven sledges crossed the Ladoga breaking ice underneath. In this case, 
strength conditions were calculated for an expendable engineering structure \cite{Kant-blokada}.

\begin{figure}[!htbp]
\centering
\includegraphics[width=0.8\textwidth]{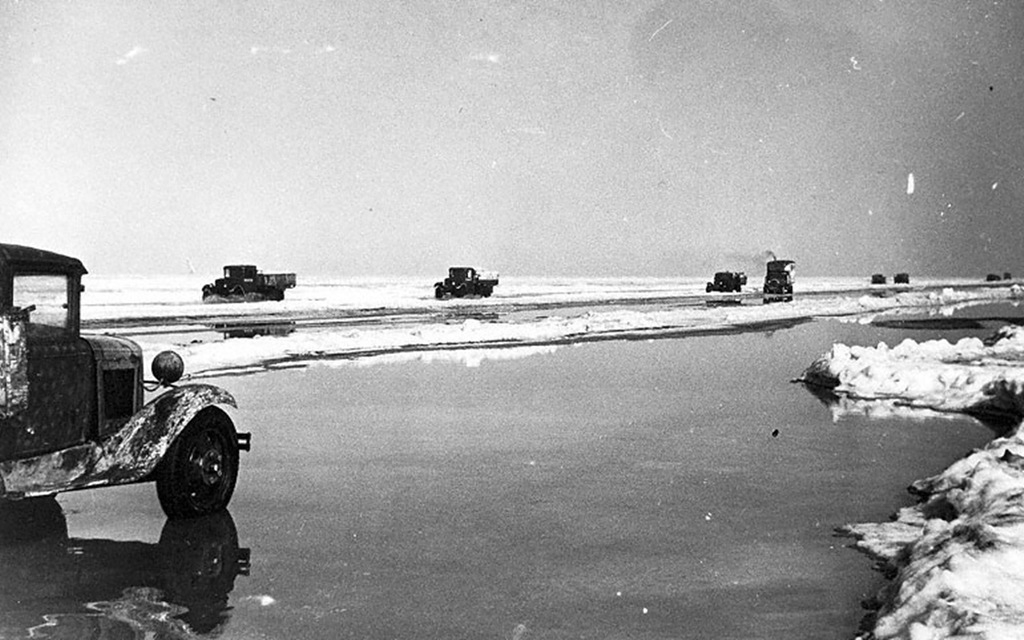}
\caption{\small The Road of Life. Lake Ladoga in spring 1942}
\end{figure}

In the beginning of 1943 the siege of Leningrad was broken, which enabled to improve Leningraders' working and living conditions a little.

The repressions on false charges persisted in Leningrad even during the siege. In winter 1942, on the charge of ``counter-revolutionary activities'' and for the ``willingness to cooperate 
with the invaders'', a big group of scientists from the city's academic institutions, University, and other higher learning institutions was arrested. Absurd charges were filed against 
them -- they allegedly colluded and planned to assassinate Stalin and form a ``\dots puppet government in Russia which would obey to Nazi'' (the role of the Prime Minister in this government 
was ascribed to V.I. Smirnov, who escaped the arrest, as he had been evacuated). Exhausting interrogations and tortures forced them to testify against their colleagues\footnote{The only two 
persons, who did not bear false witness against anyone, were V.N. Churilovsky, professor of the Leningrad Institute of Fine Mechanics and Optics and P.P. Obraztsov, a laboratory assistant 
of the University.}, the number of arrested people was growing and reached 127 men. Among the arrested scientists were V.S. Ignatovsky (shot dead), physicist, Corresponding member of the 
Academy; mathematicians, professor B.I. Izvekov and assistant professor B.D. Verzhbitsky (both died in NKVD prison). N.V. Rose, Dean of the Mathematics and Mechanics Faculty, died during 
the investigation.

Corresponding Member of the Academy N.S. Koshlyakov, head of the department of general mathematics at the University, and A.M. Zhuravsky, director of Leningrad department of the 
Mathematical Institute, were victimized on the same charges. They were sentenced to death and thereafter, this verdict was replaced by 10 years of imprisonment. A.M. Zhuravsky was 
sent to `sharazhka'\footnote{Informal name for Experimental Design Bureau, secret research and development laboratories in the Soviet Gulag labour camp system.}, where samples of 
weapons were developed. N.S. Koshlyakov stayed in a camp; in 1944, was transferred to work in Moscow at the theoretical department at Experimental 
Design Bureau. In 1951--52, they were set free, and in 1954, rehabilitated. Speaking of his involvement in defence work, A.M. Zhuravsky admitted: ``I do not work for Bolsheviks, I work for 
Russia which must be a superpower, no matter who rules it.''

\section*{IV}

Mass re-evacuation of Leningrad institutions began in the second half of 1944 and was completed by summer 1945. Some institutes, as well as certain scientists, remained to work in Moscow. 
By summer 1944, educational process began to recover in Leningrad, scientific life was resuming. Leningrad rose from the ruins. 

The war has done immense harm to the industry, urban economy, housing stock. Buildings of many institutes were damaged by artillery attacks and air strikes. Pulkovo Observatory was ruined. 
Quite a number of scientists, who returned to the city, found that they had lost their lodgings.

\begin{figure}[!htbp]
\centering
\includegraphics[width=0.8\textwidth]{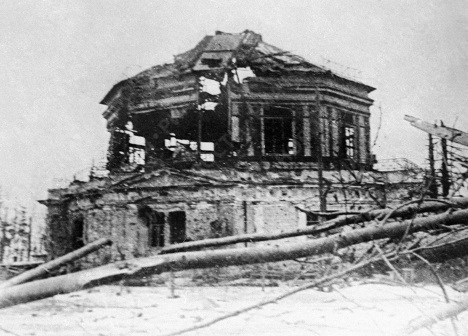}
\caption{\small Pulkovo Observatory in 1944}
\end{figure}

Many employees of institutes were killed in the lines and in the sieged city. The lack of teachers and academic specialists created such staffing bottleneck that, for example, V.I. Smirnov 
had to hold the Chair at several university departments at a time. Smirnov assigned the Chair to young professors as they were raised \cite{Smirnov}.

In the late 1940s, the practice of `bashing' scientists was resumed. They were accused of departuring Marxist positions. The best known example was the session of the All-Union Academy of 
Agricultural Sciences (1948). In the course of this session, the outstanding Soviet school of genetics was actually destroyed. Mathematicians also had to fend off ideologists' assaults.

In 1947--49, they launched a campaign to promote ``purity of physical and mathematical sciences in the socialist environment''. The constructive mathematics developed by A.A. Markov Jr. and 
cybernetics were supposed to be declared as pseudosciences. According to reminiscences of G.I. Petrashen', thanks to the presentation of A.D. Aleksandrov at an ideological seminar at the University, 
they managed to advocate constructive mathematics\footnote{Scientific community is aware of A.D. Aleksandrov's prominent services to the cause of protection of scientific biology. Thanks to his support, 
Leningrad University resumed teaching Scientific Genetics in 1950s, while in other universities in this country, including Moscow University, this happened much later.}.

\begin{figure}[!htbp]
\centering
\includegraphics[width=0.4\textwidth]{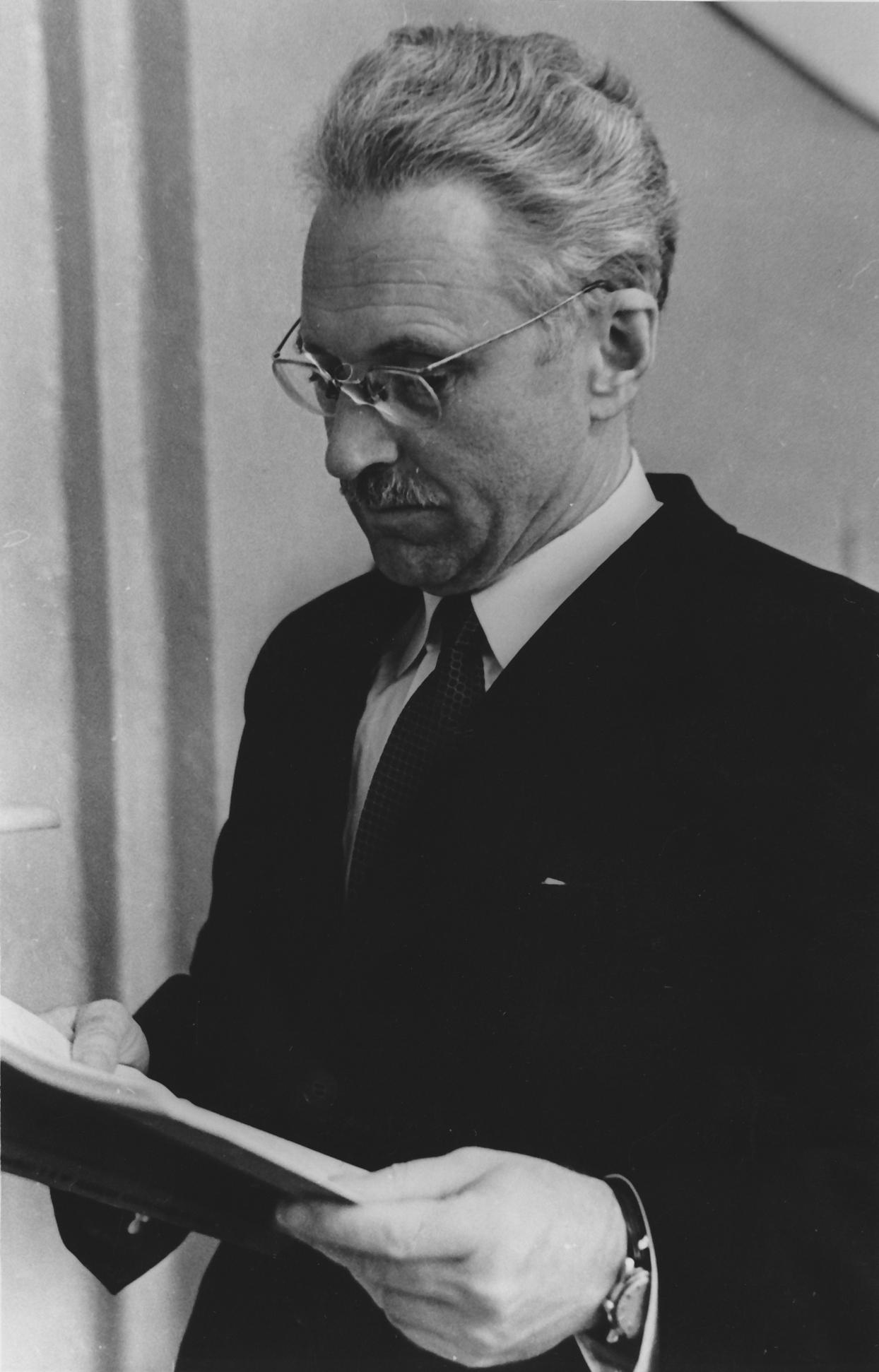}
\caption{\small A.D. Aleksandrov}
\end{figure}

In the end of 1940s -- beginning of 1950s, skies over mathematics darkened to such extent that the very existence of it as an independent science was at threat. Leading scientists felt that 
it had to be protected from anti-scientific attacks of that time. The book project of A.D. Aleksandrov, Mathematics: Its Content, Methods and Meaning, had served this purpose. It was published 
in 1953 in an edition of 350 copies with a stamp ``Printed for discussion in a limited range''. Aleksandrov himself, A.N. Kolmogorov\footnote{Andrey N. Kolmogorov (1903--1987) was a great Soviet 
mathematician, member of the Academy; one of the originators of the modern probability theory. He also achieved breakthrough results in topology, geometry, mathematical logic, classical mechanics, 
theory of turbulence, algorithmic information theory, function theory, theory of trigonometric series, measure theory, set theory, dynamical systems, functional analysis, and several other areas of 
mathematics and applications thereof.}, and M.A. Lavrentiev\footnote{Mikhail A. Lavrentiev (1900--1980) was a prominent Soviet mathematician and mechanician, founder of the Siberian Branch of the 
Academy of the USSR, member of the Academy; Vice President of the Academy of the USSR (1957--1976). He worked in the area of the complex variable theory, calculus of variations, and mathematical 
physics.} acted as editors; leningraders D.K. Faddeev, L.V. Kantorovich, V.A. Zalgaller\footnote{Victor A. Zalgaller (born in 1920) is a famous Soviet and Russian mathematician. His works are devoted 
to geometry in the large and linear programming. Since 1999, he has lives in Israel.}, O.A. Ladyzhenskaya\footnote{Olga A. Ladyzhenskaya (1922--2004) was an outstanding Soviet and Russian mathematician, 
member of the Academy. Her focus areas were partial differential equations and mathematical hydrodynamics. In her works written in co-authorship with Nina N. Uraltseva, the solution of the 19th and 
20th Gilbert problems was completed.} were among the authors of this book. In an intelligible form, the authors justified the role and 
importance of mathematics in the modern world. The goal was accomplished -- the frontal assault on mathematics was stopped\footnote{Later, the book was published in large edition in Russia (1956). 
It was also translated into English (7 editions) and Spanish. See \cite{Math-Content}.}.

The first post-war decade was marred by burst of anti-Semitism, which affected prospective students, university teachers, and staff of research institutions. Thus, in 1953, G.M. Fichtenholz was 
forced to leave the Chair at the Department of Mathematical Analysis he had himself established. To prevent an insufficiently competent person from being appointed to hold the Chair, Academician 
V.I. Smirnov took part in the competition to fill the vacancy and obtained the Chair. Thus Fichtenholz could stay at the University as professor. In 1956, Smirnov assigned the Chair to S.M. Lozinsky. 

After Stalin died (1953), no more attempts were made to destroy whole sections of mathematics, although some scientists were still plagued on account of their ethnic descent or for political reasons. 
Unspoken ethnic restrictions on admission in the University and doctorate existed until mid-1980s.

\section*{V}

The city's scientific life was reviving quite quickly after the war, as people were extremely enthusiastic. The country was in need of professionals; battle-front veterans and school leavers 
were eager to enter higher learning institutions and were very fond of their studies regardless of their stringent budget and household difficulties.

\begin{figure}[!htbp]
\centering
\includegraphics[width=0.85\textwidth]{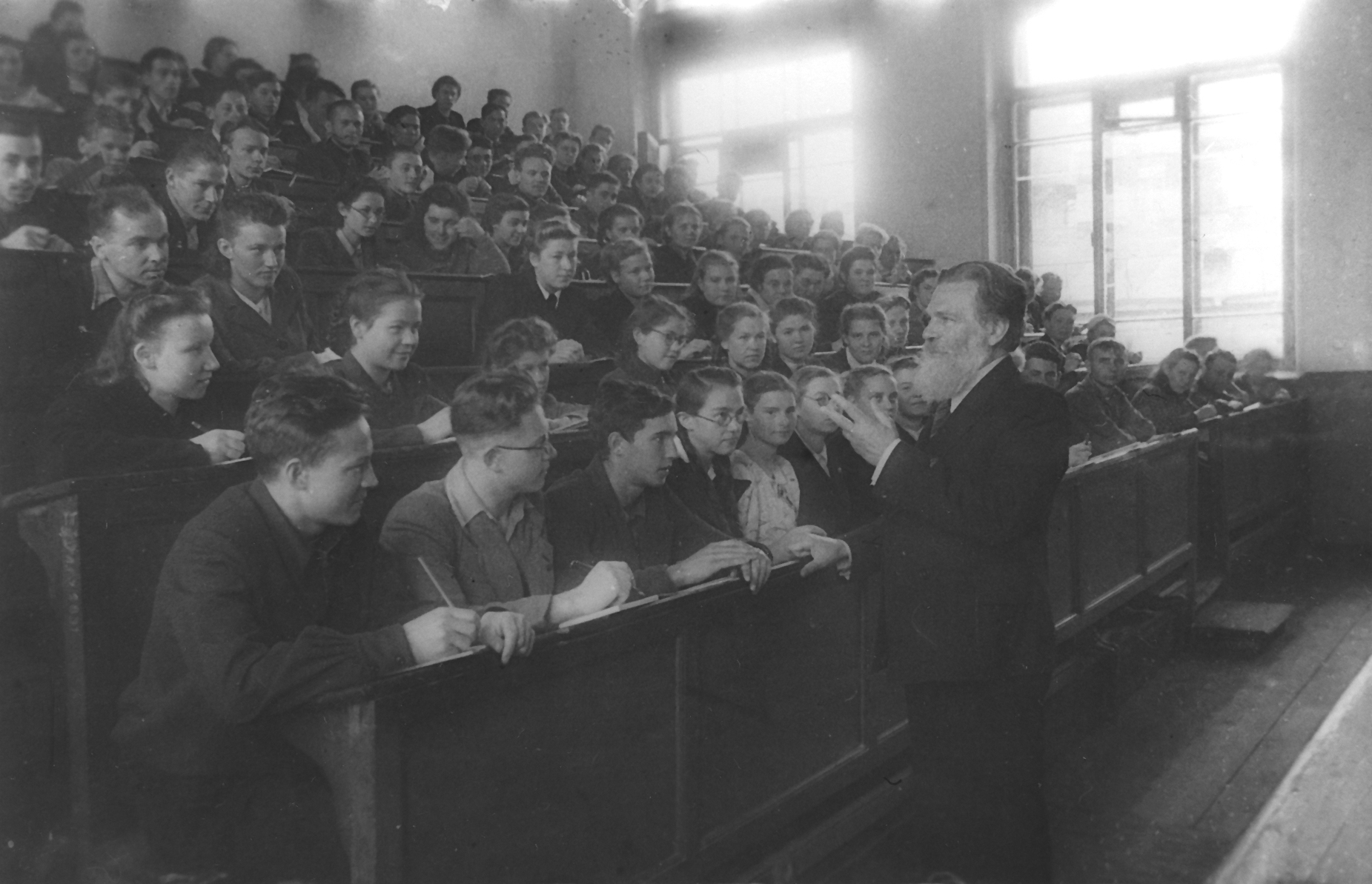}
\caption{\small September 1, 1948. First day of classes. G.M. Fichtenholz lectures to first-year students at the University}
\end{figure}

Many graduates of Leningrad University were involved in defence calculations, including the development of missile-borne nuclear weapons. In 1948, a group of 15 people under the leadership of 
L.V. Kantorovich was formed in the Leningrad Branch of the Steklov Mathematical Institute. It was charged with calculation of plutonium critical mass. They used semi-automatic calculating 
machines and later, tabulating machines to make these calculations.

The role of Leningrad mathematicians was important in other post-war scientific and technological achievements. In 1954, the first nuclear power plant was put into operation, and in 1957, 
the first artificial Earth satellite was launched.

\begin{figure}[!htbp]
\centering
\includegraphics[width=0.45\textwidth]{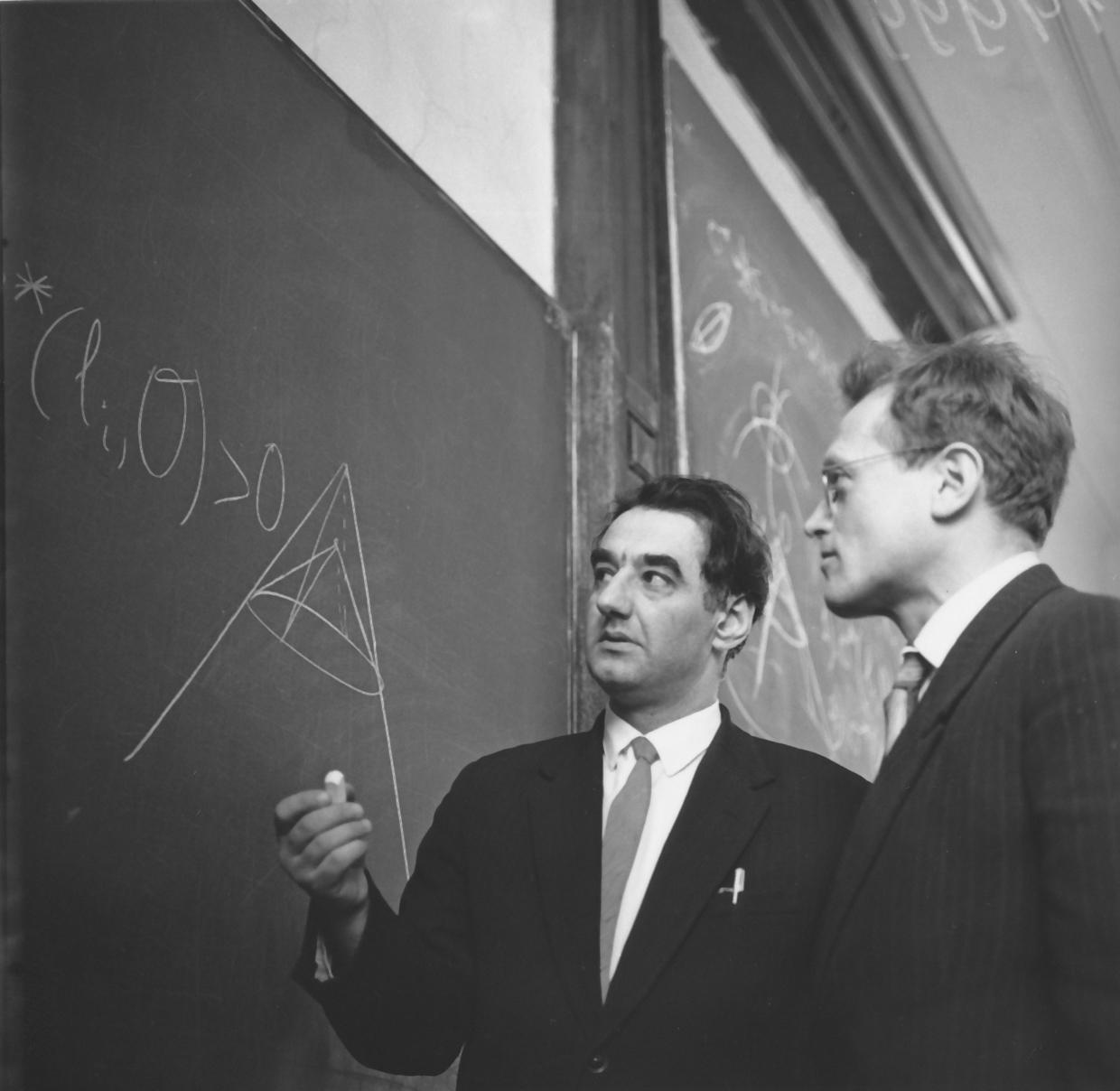}
\caption{\small V.A. Zalgaller and Yu.F. Borisov. Leningrad department of Mathematical Institute, 1963}
\end{figure}

The new generation, which joined the mathematical community after the war, gradually began to take the lead -- those were V.A. Zalgaller, O.A. Ladyzhenskaya, and thereafter, G.P. Akilov\footnote{Gleb 
P. Akilov (1921--1986) was a Soviet mathematician specialized in functional analysis. Since 1964, he worked in Novosibirsk.}, 
Z.I. Borevich\footnote{Zenon I. Borevich (1922--1995) was a Soviet and Russian mathematician specialized in algebra and number theory.}, V.A. Yakubovich\footnote{Vladimir A. Yakubovich (1926--2012) 
was a prominent Soviet and Russian mathematician, Corresponding member of Russian Academy. He made a fundamental contribution in the modern control theory and was a creator of Leningrad school of 
control theory.}, M.S. Birman\footnote{Mikhail S. Birman (1928--2009) was a prominent Soviet and Russian mathematician, founder of Leningrad school of spectral theory of operators.}, and others.

\begin{figure}[!htbp]
\centering
\includegraphics[width=0.5\textwidth]{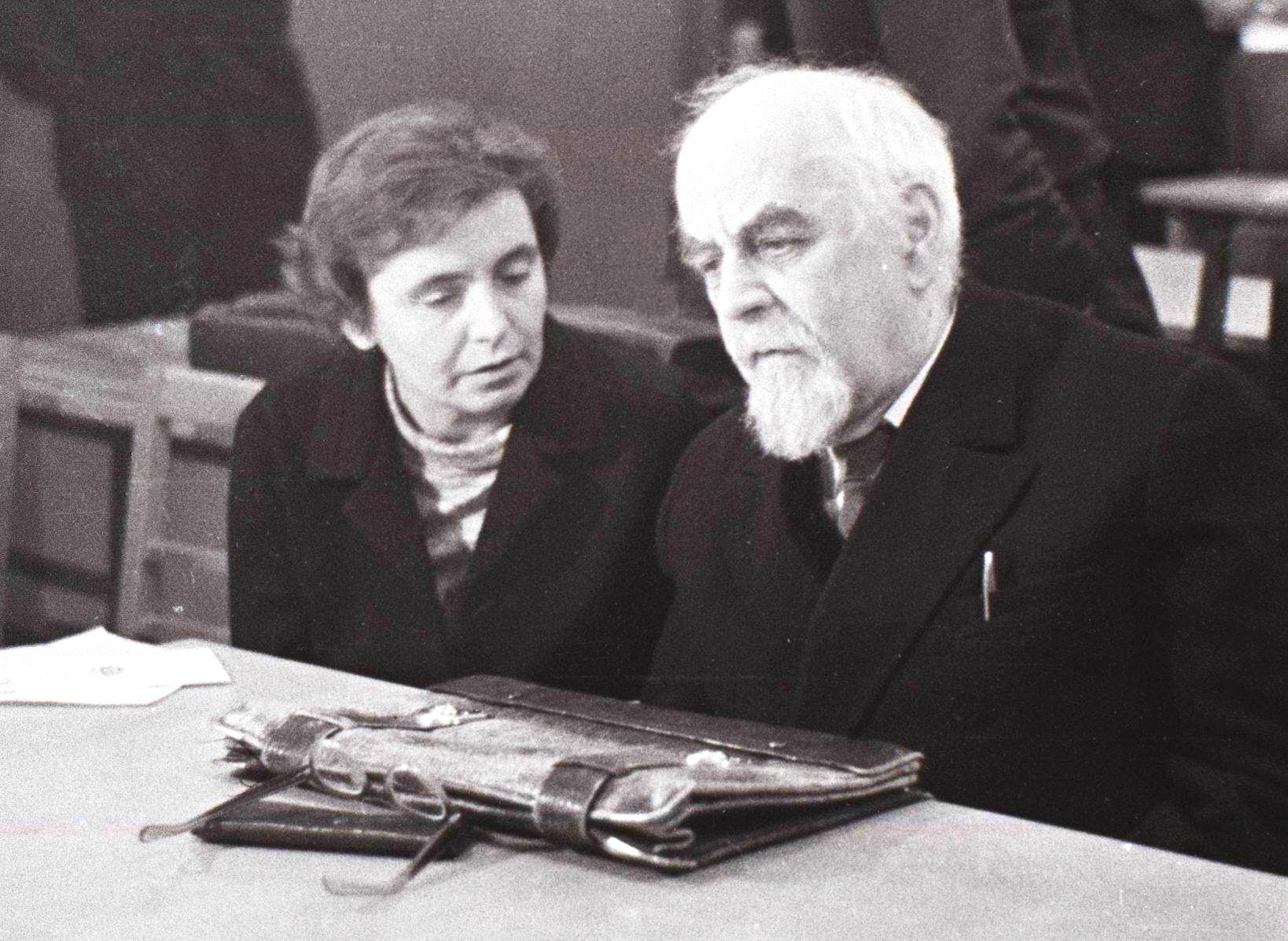}
\caption{\small O.A. Ladyzhenskaya and V.I. Smirnov, 1968}
\end{figure}

Established in 1943, the Commission for the History of Mathematical Sciences continued its work. A.N. Krylov, S.I. Vavilov, V.I. Smirnov served successively as its chairmen. The book series 
``Classics of Science'' has been published since 1945 to the present time. In 1953, the Leningrad Department of Institute of History of Natural Science and Technology was opened. In 1957, a 
Jubilee Session of the Academy of Sciences was held in Leningrad to commemorate the 250th anniversary of Euler.

\begin{figure}[!htbp]
\centering
\includegraphics[width=0.9\textwidth]{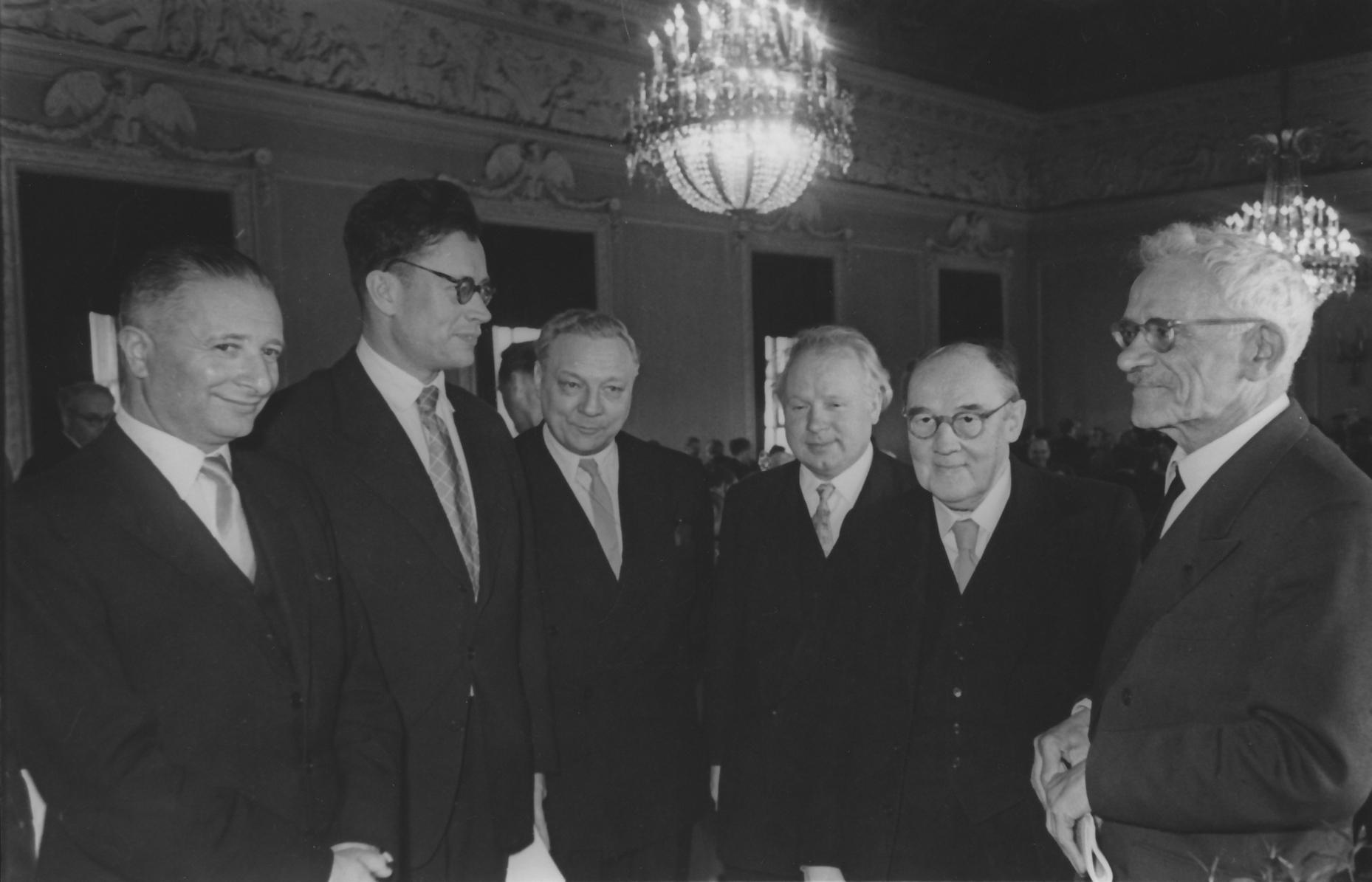}
\caption{\small Foreign participants of the Euler jubilee session. From left to right: E.~Marczewski (Poland), M.~Katetov (Czechoslovakia), H.~Grell (GDR), K.~Schr\"{o}der (GDR), W.~Sierpi\'{n}ski 
(Poland) and M.~Fr\'{e}chet (France). April 15, 1957}
\end{figure}

In the post-war period, the scope of mathematical studies was substantially expanded in Leningrad. Many conferences and scientific workshops had formative influence and were instrumental 
to attract young people to new areas of theoretical and applied mathematics. The most important research centres were in Leningrad Branch of Steklov Institute and in Leningrad State University 
as before. However, mathematical research was conducted in other higher learning institutions as well. Scientists from Leningrad were involved in training researchers and teachers for the entire 
Soviet Union as well as for foreign countries. Those were mathematicians from Leningrad who formed strong research teams in many cities of the country.

The modern Petersburg mathematical school is proud of the results in many areas: functional analysis, constructive function theory, complex variable theory, algebra and number theory, geometry 
and topology, ordinary and partial differential equations, diffraction and wave propagation theory, spectral theory, mathematical physics, dynamic systems, probability theory and mathematical 
statistics, control theory, numerical analysis, mathematical logic, computability theory and theory of computation, history of mathematics. 

L.D. Faddeev (Shao prize winner); M.L. Gromov (Abel prize winner); Fields medalists G.Ya. Perelman and S.K. Smirnov; and many other outstanding mathematicians are among those who were 
trained by this school.


\begin{thebibliography}{8}
\bibitem{MatPet} 
Mathematical Petersburg. History, science, sights / G.I. Sinkevich, editor-compiler; A.I. Nazarov, scientific editor. SPb: Educational projects. 2018. 336 pp. In Russian.

\bibitem{Krylov}
Academic Science in St. Petersburg in the 18th-20th centuries. Historical Essays /Editor-in-Charge J.I. Alferov; St. Petersburg Branch of S.I. Vavilov Institute of History of Natural Science and 
Technology. SPb: Nauka. 2003. 605 pp. In Russian.

\bibitem{Luzin}
The Case of Academician Nikolai Nikolaevich Luzin / Ed. by Sergei S. Demidov, Boris V. L\"{e}vshin. AMS: History of Mathematics. V. 43, 2016. 416 pp.

\bibitem{LenMatFront}
At the Leningrad Mathematical Front. M.-L.: GSEI. 1931. 44 pp. In Russian.

\bibitem{blokada}
Ezhov V.A., Mavrodin V.V. Leningrad University During the Great Patriotic War. L.: Leningrad State University Publishers. 1975. 88 pp. In Russian.

\bibitem{Kant-blokada}
Petrenko I.V. Leonid Vital’evich Kantorovich // As Recollected by Leningraders. In Russian. Available at http://www.liveinternet.ru/users/3652449/post135225542/ 

\bibitem{Smirnov}
Apushkinskaya D.E., Nazarov A.I. Vladimir Ivanovich Smirnov (1887--1974) // Complex Variables and Elliptic Equations. 2018. V. 63. 897--906.

\bibitem{Math-Content}
Aleksandrov A.D., Kolmogorov A.N., Lavrentiev M.A., ed: Mathematics: Its Content, Methods and Meaning (3 Volumes in One). Courier Corporation, 2012. 1120 pp.

\end{thebibliography}
\end{document}